\documentclass[12pt, english]{article}
\usepackage[T1]{fontenc}
\usepackage{amsmath}
\usepackage{amsthm}
\usepackage{amssymb}
\usepackage{enumitem}
\usepackage{xcolor}
\usepackage{tikz-cd}
\usepackage[blocks]{authblk}
\usepackage[margin=1in,includefoot,headheight=0pt,headsep=0pt,footskip=30pt]{geometry}
\usepackage{babel}
\usepackage{newtxtext}
\usepackage{mathcommand}
\usepackage{fixdif}
\usepackage{physics2}
\usepackage[%
pagebackref, 
breaklinks, 
colorlinks=true, linkcolor=blue, citecolor=blue, urlcolor=blue
]{hyperref}
\usephysicsmodule{ab}

\numberwithin{equation}{section}
\allowdisplaybreaks

\DeclareDocumentMathCommand{\odv}{m}{\frac{\d}{\d{#1}}}
\DeclareDocumentMathCommand{\dbrack}{m}{\mathopen{}[\mkern-4mu[#1]\mkern-4mu]\mathclose{}}
\DeclareDocumentMathCommand{\dparen}{m}{\mathopen{}(\mkern-4mu(#1)\mkern-4mu)\mathclose{}}
\declaremathcommand{\C}{{\mathbb{C}}}
\declaremathcommand{\Z}{{\mathbb{Z}}}
\declaremathcommand{\R}{{\mathbb{R}}}
\declaremathcommand{\N}{{\mathbb{N}}}
\declaremathcommand{\A}{{\mathcal{A}}}
\declaremathcommand{\Y}{{\mathcal{Y}}}

\DeclareMathOperator{\End}{End}
\DeclareMathOperator{\wt}{wt}
\DeclareMathOperator{\Hom}{Hom}
\DeclareMathOperator{\id}{id}
\DeclareMathOperator{\opp}{op}
\DeclareMathOperator{\lr}{lr}
\DeclareMathOperator{\Nat}{Nat}
\DeclareMathOperator{\Dinat}{Dinat}
\DeclareMathOperator{\Span}{Span}
\DeclareMathOperator{\Image}{Im}
\DeclareMathOperator{\Ker}{Ker}
\DeclareMathOperator{\Rex}{Rex}
\DeclareMathOperator{\Brex}{Brex}
\DeclareMathOperator{\Coend}{Coend}
\DeclareMathOperator{\Log}{Log}
\DeclareMathOperator{\Arg}{Arg}
\DeclareMathOperator{\ev}{ev}

\theoremstyle{plain}

\newtheorem{mainthm}{Theorem}

\newtheorem{lemma}{Lemma}[section]
\newtheorem{theorem}[lemma]{Theorem}
\newtheorem{corollary}[lemma]{Corollary}
\newtheorem{proposition}[lemma]{Proposition}

\theoremstyle{definition}
\newtheorem{definition}[lemma]{Definition}

\newtheorem{example}[lemma]{Example}

\theoremstyle{remark}

\newtheorem{remark}[lemma]{Remark}

\title{Representability of Intertwining Operators for Weak Modules of Vertex Operator Algebras}

\author{Xu Gao \\ Email: gaoxu@shmtu.edu.cn }
\affil{Department of Mathematics, School of Sciences, Shanghai Maritime University,
	1550 Haigang Avenue, Pudong New Area,
	Shanghai 201306, China}

\author{Chao Yang   \\ Email: yangcabc@163.com }
\affil{ School of Mathematics,  Southwest Jiaotong University,
	Chengdu 611756, China   }

\date{}

\begin{document}
	\maketitle
	
	\begin{abstract}
		
		Let $V$ be a vertex operator algebra. For every weak $V$-module $W$, we construct a weak $V \otimes V$-module $\mathcal{R}(W)$ representing the bifunctor $\mathcal{I}\binom{W}{-\;-}$ of intertwining operator spaces, in the sense that there exists a natural isomorphism 
		\[
		\mathcal{I}\binom{W}{-\;-} \cong \Hom_{V\otimes V} \ab\big(-\otimes -, \; \mathcal{R}(W)).
		\]
		As applications of this natural isomorphism, we first establish a Yoneda-type lemma for spaces of intertwining operators, and then identify the unique largest locally restricted submodule $\mathcal{R}_{\lr}(V^\prime)$ of $\mathcal{R}(V^\prime)$ with the coend of the bifunctor $-\otimes-^\prime$:
		\[
		\mathcal{R}_{\lr}(V^\prime) \cong \int^{M} M \otimes M^\prime,
		\]
		where $M$ ranges over restricted $V$-modules.
		Finally, the connection between $\mathcal{R}(W^\prime)$ and the regular representation $\mathfrak{D}_{\xi}(W)$ constructed in \cite{Li02} is established, where $\xi \in \C^{\times}$. In particular, if $V$ is a $C_2$-cofinite vertex operator algebra and $W$ is a restricted $V$-module, then $\mathcal{R}(W^\prime)$ is also a restricted $V \otimes V$-module and is canonically isomorphic to $\mathfrak{D}_{\xi}(W)$.

	\end{abstract}
	
	\tableofcontents
	
	\section{Introduction}
	
	Let $V$ be a vertex operator algebra. For weak $V$-modules $M, N, W$, let $\mathcal{I}\binom{W}{M\;N}$ denote the space of \emph{logarithmic intertwining operators} of the indicated type. This gives rise to a bifunctor, contravariant in both variables,
	\[
	\mathcal{I}^{W}(-, -) := \mathcal{I}\binom{W}{-\;-}
	\colon
	(\mathfrak{M}^w_V)^{\opp}\times
	(\mathfrak{M}^w_V)^{\opp}
	\longrightarrow \mathbf{Vect},
	\]
	where $\mathfrak{M}^w_V$ is the category of weak $V$-modules and $\mathbf{Vect}$ is the category of complex vector spaces.
	It is then natural to ask whether this bifunctor can be \emph{represented} by a weak $V\otimes V$-module $\mathcal{R}(W)$. That is, whether there exists a natural isomorphism of bifunctors:
	\begin{equation}\label{eq:intertwiner-hom-intro}
		\mathcal{I}^{W}(-, -) \cong 
		\Hom_{V \otimes V} \ab\big(-\otimes-, \mathcal{R}(W)).
	\end{equation}
	This can be viewed as an external version of the usual representability problem. 
	In analogy with Yoneda notation, we write $\mathcal{H}_X=\Hom_{V\otimes V}(-\otimes-,X)$ for a weak $V\otimes V$-module $X$.

	In the category of ordinary $V$-modules, this has been answered by Li's \emph{regular representations}. In \cite{Li02}, for each weak $V$-module $W$ and $\xi \in \C^{\times}$, a weak $V\otimes V$-module $\mathfrak{D}_{\xi}(W)$ is constructed from the full dual space $W^*$, and it is shown that there exists an isomorphism
	\begin{equation}\label{eq:li-hom-intro}
		\mathcal{I}\binom{W^\prime}{M \; N}
		\cong
		\Hom_{V\otimes V} \ab\big(M\otimes N, \mathfrak{D}_{\xi}(W))
	\end{equation}
	for all ordinary $V$-modules $M, N, W$.
	This construction has led to a series of applications.
	For example, it is used in \cite{Li01} to simplify the proof of the fusion rule theorem from \cite{FZ92,Li99a}; in \cite{Li01,Li23}, the construction of $\mathfrak{D}_\xi(W)$ has been connected to the theory of generalized Zhu algebras in \cite{Zhu96,DLM98,DJ08}.
	Twisted analogues of this construction have been studied in \cite{LS23,Zy25}.
	
	However, the original proof of the isomorphism \eqref{eq:li-hom-intro} relies on the construction of intertwining operators from \emph{$P(\xi)$-intertwining maps}. It is therefore difficult to establish such a description for arbitrary weak modules. In contrast, the initial goal of this paper is to construct the representing object $\mathcal{R}(W)$ for arbitrary weak $V$-module $W$ directly from intertwining operators.
	
	The idea is simple: we package all intertwining operators of type $\binom{W}{M\;N}$ for all weak $V$-modules $M, N$ into a single $\mathcal{R}(W)$:
	\[
	\mathcal{R}(W)
	=
	\Span_{\C}
	\ab\{
	\mathcal{Y}(a, y) b
	\; \middle| \;
	M, N \in \mathfrak{M}^w_V, \
	a \in M, \ b \in N, \
	\mathcal{Y} \in \mathcal{I}\binom{W}{M\,N}
	\}
	\subset W[\log y] \{y\}.
	\]
	Then we define a weak $V\otimes V$-module structure on it
	\[
	Y_{\mathcal{R}(W)}(-, x) \colon \; (V \otimes V) \otimes \mathcal{R}(W) \longrightarrow \mathcal{R}(W) \dparen{x}
	\]
	by
	\[
	u \otimes v \otimes \Y(a, y) b \longmapsto \Y \ab\big(Y_{M}(u, x) a, y) Y_{N}(v, x) b
	\]
	for all $u, v \in V$ and $\Y(a, y) b \in \mathcal{R}(W)$. 
	Our first main result is the following.
	\begin{mainthm} \label{thm:hom-realization-intro}
		Let $V$ be a vertex operator algebra, and let $W$ be a weak $V$-module.
		Then the pair $(\mathcal{R}(W), Y_{\mathcal{R}(W)})$ forms a weak $V \otimes V$-module.
		For any weak $V$-modules $M, N$, the map
		\begin{equation} \label{eq:hom-realization-map-intro}
			\psi^{W}_{M, N}\colon \mathcal{I}^{W}(M, N) \longrightarrow \Hom_{V \otimes V} \ab\big(M \otimes N, \; \mathcal{R}(W)), \quad \Y \longrightarrow \ab\big(a\otimes b \mapsto \Y(a, y) b)
		\end{equation}
		is a linear isomorphism.
		Moreover, these isomorphisms are natural in $M$ and $N$. 
		In other words, the bifunctor $\mathcal{I}^{W}(-, -)$ is isomorphic to $\mathcal{H}_{\mathcal{R}(W)}$.
	\end{mainthm}
	
	It is then natural to ask how $\mathcal{R}(W^\prime)$ is related to Li's regular representation $\mathfrak{D}_{\xi}(W)$. This question is not easy to answer in general. However, we can still establish a canonical isomorphism between the largest bi-generalized $V\otimes V$-submodule of $\mathcal{R}(W^\prime)$ and that of $\mathfrak{D}_{\xi}(W)$ when $W$ is a \emph{restricted} $V$-module (see Theorem~\ref{thm:li-regular-comparison}). 
	A more satisfactory answer is available in the $C_2$-cofinite case: If $V$ is a $C_2$-cofinite vertex operator algebra and $W$ is a restricted $V$-module, then $\mathcal{R}(W^\prime)$ is also a restricted $V\otimes V$-module and is canonically isomorphic to $\mathfrak{D}_{\xi}(W)$.
	Another natural question about $\mathcal{R}(W)$ is whether it can be constructed without using intertwining operators. More precisely, we ask which subspace of $W[\log y] \{y\}$ coincides with $\mathcal{R}(W)$.
	We will address this question in our forthcoming paper.

	To obtain applications of Theorem~\ref{thm:hom-realization-intro}, we establish the following generalized version of the Yoneda lemma.
	
	\begin{mainthm}\label{thm:external-yoneda-intro}
		The assignment $X \mapsto \mathcal{H}_X$ defines an embedding from the category $\mathfrak{M}^w_{V \otimes V}$ into the category of bifunctors $(\mathfrak{M}^w_V)^{\opp}\times(\mathfrak{M}^w_V)^{\opp} \to \mathbf{Vect}$. 
		In particular, for any weak $V\otimes V$-modules $X, Z$, composition induces an isomorphism
		\begin{equation}\label{eq:external-yoneda-intro}
			\Hom_{V\otimes V} \ab\big(X, Z)
			\cong \Nat \ab\big(\mathcal{H}_X, \mathcal{H}_Z).
		\end{equation}
		Moreover, if $X$ is locally restricted, then
		\begin{equation}\label{eq:restricted-external-yoneda-intro}
			\Hom_{V\otimes V} \ab\big(X, Z)
			\cong \Nat_{(\mathfrak{M}^r_V)^{\opp}\times(\mathfrak{M}^r_V)^{\opp}} \ab\big(\mathcal{H}_X, \mathcal{H}_Z).
		\end{equation}
	\end{mainthm}

	Here and throughout this paper, we write $\Nat_{\mathcal{C}}(F,G)$ for natural transformations between $F|_{\mathcal{C}}$ and $G|_{\mathcal{C}}$, omitting the subscript when the domain is clear.
	
	Theorems~\ref{thm:hom-realization-intro} and~\ref{thm:external-yoneda-intro} form the core technical tools for our applications.
	As a first consequence, we obtain the following Yoneda-type correspondence for intertwiner bifunctors.

	\begin{mainthm} \label{Yoneda-lem-intro}
		Let $V$ be a vertex operator algebra, and let $W, M$ be weak $V$-modules. Then there is a linear isomorphism
		\begin{equation*}
			\operatorname{Nat}\bigl( \mathcal{I}^W, \mathcal I^M \bigr) \cong \operatorname{Hom}_{V \otimes V}\bigl( \mathcal R(W), \mathcal R(M) \bigr).
		\end{equation*}
		As a consequence, the functors $\mathcal I^W$ and $\mathcal I^M$ are isomorphic if and only if $W\cong M$ as weak $V$-modules.
	\end{mainthm}

	\begin{center}
		\textbf{\textasteriskcentered\qquad\textasteriskcentered\qquad\textasteriskcentered}
	\end{center}
	
	When $V$ is \emph{regular}, a Peter-Weyl-type decomposition of $\mathfrak{D}_\xi(V)$ is established in \cite{Li02,Li02b}. More generally, under the $C_2$-cofiniteness assumption, $\mathfrak{D}_\xi(V)$ is interpreted in terms of conformal blocks and is identified with the \emph{coend} of the bifunctor
	\[
	-\otimes -' \colon \mathfrak{M}^r_V \times (\mathfrak{M}^r_V)^{\opp} \to \mathfrak{M}^r_{V \otimes V},
	\]
	where $\mathfrak{M}^r_V$ is the category of restricted $V$-modules and $-'$ is the restricted dual functor; see \cite{GZ25,GZ26}. This coend realization of $\mathfrak{D}_\xi(V)$ is a categorical generalization of the Peter-Weyl-type decomposition. We point out that the coend realizations in \cite{GZ25,GZ26} rely on the fact that the category of restricted modules over a $C_2$-cofinite vertex operator algebra is a finite linear abelian category.
	Concretely, this finiteness property is used to guarantee properties such as the existence of Deligne tensor products and of projective generators.

	In general, $\mathfrak{D}_\xi(V)$ or $\mathcal{R}(V')$ need not be a restricted $V\otimes V$-module, and the index category $\mathfrak{M}^r_V$ is not necessarily a finite abelian category. 
	Consequently, the argument from \cite{GZ25,GZ26} does not apply. 
	However, as an application of Theorem~\ref{thm:hom-realization-intro} and Theorem~\ref{thm:external-yoneda-intro}, a coend formula can still be obtained.

	We now explain the identification of $\mathcal{R}(V')$ with the coend of the bifunctor $ -\otimes -'$.
	By the definition of a coend (see Section~\ref{Section-coend}), to establish this identification we first need a family of morphisms
	\[
	\pi_M\colon M\otimes M' \to \mathcal R(V^\prime), \quad M \in \mathfrak{M}^r_V,
	\]
	such that for every morphism $f\colon M \to N$ in $\mathfrak{M}^r_V$,  the following diagram commutes:
	\[
	\begin{tikzcd}[row sep=small]
		& M\otimes M' \arrow[dr, "\pi_M"] & \\
		M\otimes N'
		\arrow[ur, "M\otimes f'"]
		\arrow[dr, "f\otimes{N'}"']
		& & \mathcal R(V^\prime) \\
		& N\otimes N' \arrow[ur, "\pi_N"'] &
	\end{tikzcd}
	\]
	Such a family is called a dinatural transformation from the bifunctor $-\otimes -'$ to the object $R(V^\prime)$. 
	Second, we need to verify that this dinatural transformation is initial among all such transformations. 
	More precisely, for any restricted $V\otimes V$-module $X$ and any dinatural transformation 
	\[
	\eta_M\colon M\otimes M' \to X, \quad M \in \mathfrak{M}^r_V,
	\]
	there is a unique module homomorphism $ R(V^\prime) \to X$ making the triangle commute for all $M\in\mathfrak{M}^r_V$:
	\[
	\begin{tikzcd}[column sep=large]
		M\otimes M'
		\arrow[r, "\pi_M"]
		\arrow[dr, "\eta_M"']
		&  R(V^\prime) \arrow[d, dashed, "{\exists!}"] \\
		& X
	\end{tikzcd}
	\]
	It is often to denote the coend of a bifunctor in integral notation. In this case, it is written as
	\[
	R(V^\prime) = \int^{M\in\mathfrak{M}^r_V} M\otimes M'.
	\]

	By  Theorem~\ref{thm:hom-realization-intro} and the symmetries of intertwining operators, we have
	\[
	\Hom_V(M, N)
	\cong \Hom_{V\otimes V}
	\ab\big(M\otimes N', \mathcal{R}(V')),
	\]
	for all restricted $V$-modules $M, N$ (see Lemma~\ref{lem:hom-kernel-realization}).
	Let $\mathcal{R}_{\lr}(V')$ be the largest locally restricted submodule of $\mathcal{R}(V')$.
	Then any homomorphism $M\otimes N' \to \mathcal{R}(V')$ factors through it.
	Consequently, we get a natural isomorphism:
	\[
	\Hom_V(-, -) \cong \Hom_{V\otimes V} \ab\big( -\otimes -', \; \mathcal{R}_{\lr}(V')).
	\]
	Combining this with the natural isomorphism \eqref{eq:restricted-external-yoneda-intro}, we obtain, for any weak $V\otimes V$-module $X$, the following isomorphisms:
	\[
	\begin{aligned}
		\Hom_{V\otimes V} \ab\big(\mathcal{R}_{\lr}(V'), X)
		&\cong
		\Nat_{(\mathfrak{M}^r_V)^{\opp}\times(\mathfrak{M}^r_V)^{\opp}} \ab\big(\mathcal{H}_{\mathcal{R}_{\lr}(V')}, \mathcal{H}_X)
		\\
		&\cong
		\Nat_{(\mathfrak{M}^r_V)^{\opp}\times\mathfrak{M}^r_V} \ab\big(\Hom_V(-,-), \mathcal{H}_X(-, -^\prime)),
	\end{aligned}
	\]
	where the last is obtained by \eqref{eq:intertwiner-hom-intro} and replacing the second input by its restricted dual.
	On the other hand, using Corollary~\ref{cor:natural-dinatural-bijection}, we obtain
	\begin{equation} \label{Intro-Dinat=Hom}
		\operatorname{Dinat}( -\otimes-^\prime, X) \cong \Hom_{V\otimes V} \ab\big(\mathcal{R}_{\lr}(V'), X),
	\end{equation}
	where $\operatorname{Dinat}(-\otimes-^\prime, X)$ denotes the space of dinatural transformations from $-\otimes -^\prime$ to $X$.
	This means that every dinatural transformation from the bifunctor $-\otimes -'$ to the object $X$   
	is induced by a unique $V\otimes V$-module homomorphism from $\mathcal {R}(V^\prime)$ to $X$.
	This gives precisely the universal property of the coend, hence
	\begin{mainthm}\label{thm:restricted-coend-intro}
		Let $V$ be a vertex operator algebra. Then  $\displaystyle \mathcal{R}_{\lr}(V') \cong \int^{M\in\mathfrak{M}^r_V} M\otimes M'$.
	\end{mainthm}

	\begin{center}
		\textbf{\textasteriskcentered\qquad\textasteriskcentered\qquad\textasteriskcentered}
	\end{center}

	This paper is organized as follows. Section~\ref{sec:preliminaries} recalls fundamental concepts related to modules over vertex operator algebras and intertwining operators, introduces the trifunctor $\mathcal{I} \binom{-}{- \; -}$, and reviews the properties of $V \otimes V$-modules. In Section~\ref{sec:regular-construction}, for every weak $V$-module $W$, we construct a weak $V \otimes V$-module $\mathcal{R}(W)$ and prove Theorem~\ref{thm:hom-realization-intro}. In Section~\ref{sec:c2-cofinite-case}, we work in the $C_2$-cofinite setting. We first review properties of the weak module category of a $C_2$-cofinite vertex operator algebra, and then prove the finiteness of $\mathcal{R}(W)$ for restricted $V$-module $W$. Finally, using the Deligne tensor product method for finite linear abelian categories, we prove a Yoneda-type lemma for intertwining operator spaces. In Section~\ref{sec:yoneda-and-coend}, we prove Theorems~\ref{thm:external-yoneda-intro} and~\ref{thm:restricted-coend-intro} in a more general setting. Our main tool is the external Yoneda lemma (Lemma~\ref{lem:external-tensor-density}).
	Section~\ref{sec:li-comparison} first recalls the construction of the regular representation $\mathfrak{D}_{\xi}(W)$ from \cite{Li02}, and then establishes a relation between $\mathcal{R}(W^\prime)$ and $\mathfrak{D}_\xi(W)$ for a restricted $V$-module $W$.

	\subsection*{Conventions}
	
	All vector spaces, linear maps, and tensor products are over $\C$.
	We write $\mathbf{Vect}$ for the category of complex vector spaces and linear maps. 
	All categories and functors are $\C$-linear, namely enriched over $\mathbf{Vect}$. In particular, by an abelian category we mean a $\C$-linear abelian category. 
	By a \textbf{presheaf} we mean a contravariant functor with values in $\mathbf{Vect}$.
	All finite-dimensional complex vector spaces form a full abelian subcategory of $\mathbf{Vect}$, which we denote by $\mathbf{Vect}_f$.
	Except for $\xi$, which is a nonzero complex number, all variables are formal variables, and their logarithms are formal symbols.
	
	\subsection*{Notation}
	
	\begin{itemize}
		\item $M, N, M_i, N_i$; $W, W_i$; $X, Z, X_i$: Input modules; target modules; modules over a tensor product of two VOAs, respectively.
		
		\item $\mathfrak{M}^w_V$, $\mathfrak{M}^g_V$, $\mathfrak{M}^r_V$: Categories of weak, generalized, and restricted $V$-modules, respectively (see Section~\ref{sec:preliminaries}).
		
		\item $M^*$, $M^\prime$, $f^\prime$: Full dual, restricted dual, and restricted transpose, respectively (see Section~\ref{sec:preliminaries}).
		
		\item $\mathcal{I}\binom{W}{M\;N}$, $\mathcal{I}^{W}(M, N)$: Space of logarithmic intertwining operators with inputs $M$, $N$ and target $W$ (see Section~\ref{sec:preliminaries}).
		
		\item $\mathcal{R}(W)$: Weak $V\otimes V$-module constructed from intertwining operators with target $W$ (see \eqref{eq:regular-construction}).
		
		\item $\psi^{W}_{M, N}$: Natural isomorphism from intertwining operators to $V\otimes V$-module homomorphisms (see \eqref{eq:hom-realization-map}).
		
		\item $h_X$, $\mathcal{H}_X$: The Yoneda presheaf $\Hom(-,X)$ and external Hom bifunctor $\Hom(-\otimes-,X)$, with domains specified where used.
		
		\item $\Theta$: Canonical map $\beta \mapsto \mathcal{H}(\beta)$ from module homomorphisms to natural transformations (see Lemma~\ref{lem:external-tensor-density}).
		
		\item $\mathcal{A}\boxtimes\mathcal{B}$: Deligne tensor product of linear abelian categories (see Section~\ref{sec:c2-cofinite-case}).
		
		\item $\mathcal{A}\odot\mathcal{B}$: Class of external tensor products $M\otimes N$ with $M\in\mathcal{A}$ and $N\in\mathcal{B}$ (see Section~\ref{subsec:generated-modules}).
		
		\item $\mathcal{G}_{\mathcal{A}, \mathcal{B}}$: Category of weak $V\otimes U$-modules generated by $\mathcal{A}\odot\mathcal{B}$ (see Section~\ref{subsec:generated-modules}).
		
		\item $\mathcal{S}_{\mathcal{A}, \mathcal{B}}(X)$: Largest submodule of $X$ generated by $\mathcal{A}\odot\mathcal{B}$ (see Section~\ref{subsec:generated-modules}).
		
		\item $\mathcal{R}_{\mathcal{A}, \mathcal{B}}(W)$: The submodule $\mathcal{S}_{\mathcal{A}, \mathcal{B}}\ab\big(\mathcal{R}(W))$, for $U=V$ (see Section~\ref{subsec:generated-modules}).
		
		\item $\mathcal{R}_{\lr}(V^\prime)$: Largest locally restricted submodule of $\mathcal{R}(V^\prime)$ (see the end of Section~\ref{sec:yoneda-and-coend}).
		
		\item $\mathfrak{D}_\xi(W)$: Li's regular representation inside the full dual $W^*$, for $\xi\ne0$ (see \eqref{eq:li-functional-condition}).
	\end{itemize}
	
	\section{Preliminaries}\label{sec:preliminaries}
	
	\subsection{Modules and contragredients}
	
	We adopt the standard framework for vertex operator algebras as presented in \cite{FLM88, LL04},
	and denote a vertex operator algebra by $(V, Y(-, x), \mathbf{1}, \omega)$.
	By a \textbf{weak $V$-module} we mean a vector space $M$ equipped with a linear map
	\[
	Y_M \colon V \otimes M \longrightarrow M\dparen{x}, \quad v \otimes m \longrightarrow Y_M(v, x) m,
	\]
	satisfying $Y_M(\mathbf{1}, x) = \id_M$, and the Jacobi identity.
	In this paper, we write
	\[
	Y_M(\omega, x) = \sum_{n \in \Z} L_M(n) x^{-n-2}
	\]
	and, when no confusion arises, we write $L(n)$ for $L_M(n)$.
	
	Let $V$ be a vertex operator algebra, and let $M$ be a weak $V$-module.
	For any $\lambda \in \C$, let $M_{(\lambda)}$ denote the generalized $\lambda$-eigenspace of $L(0)$ on $M$.
	When $M_{(\lambda)}$ is nonzero, $M_{(\lambda)}$ is called a \textbf{generalized weight space} of $M$, and $\lambda$ is called a \textbf{generalized weight} of $M$.
	
	\begin{definition}
		
		Let $V$ be a vertex operator algebra.
		
		\begin{enumerate}[label=(\arabic*),ref=\arabic*]
			\item A \textbf{ generalized $V$-module} is a weak $V$-module $M$ equipped with the generalized weight space decomposition $M = \bigoplus_{\lambda \in \C} M_{(\lambda)}$.
			
			\item A \textbf{grading-restricted generalized $V$-module} (or simply a \textbf{restricted $V$-module}) is a generalized $V$-module $M$ such that each generalized weight space $M_{(\lambda)}$ is finite-dimensional, and for every $\lambda \in \C$, $M_{(\lambda-n)} = 0$ for all sufficiently large positive integers $n$.
			
			\item  An \textbf{ordinary $V$-module} is a restricted $V$-module on which $L(0)$ acts semisimply.
			
			\item Let $ \mathfrak{M}^w_V, \; \mathfrak{M}^g_V$ and $ \mathfrak{M}^r_V $  denote the categories of weak, generalized, and restricted $V$-modules, respectively.
			
		\end{enumerate}
		
	\end{definition}
	
	\begin{definition}
		Let $V$ be a vertex operator algebra. A weak $V$-module $M$ is said to be $C_1$-cofinite if the quotient space $M/C_1(M)$ is finite-dimensional,
		where
		\[
		C_1(M) = \operatorname{span}\ab\{ v_{-1}m \ | \ v \in V, \  \wt (v) > 0, \ m \in M \}.
		\]
	\end{definition}
	
	\noindent \textbf{Contragredient modules}.
	Let $V$ be a vertex operator algebra, and let $(M, Y_M)$ be a  weak  $V$-module.
	Denote by $M^*$ the full dual space of $M$. We define a linear map
	\[
	Y_{M^*} \colon V \otimes M^* \to M^*\dbrack{x, x^{-1}}
	\]
	by the identity
	\begin{equation}\label{eq:full-dual-action}
		\ab\big\langle Y_{M^*}(u, x) f, \; m\rangle = \ab\big\langle f, \; Y_M^{\circ}(u, x) m\rangle, \quad \forall \; u \in V, m \in M, f \in M^*.
	\end{equation}
	Here and below,
	\begin{equation}
		Y_{M}^{\circ}(u, x) = Y_{M}(e^{xL(1)}(-x^{-2})^{L(0)}u, x^{-1}).
	\end{equation}
	In general, the pair $(M^*, Y_{M^*})$ is not a weak $V$-module.
	
	Now suppose that $(M, Y_M)$ is a restricted $V$-module.
	The restricted dual space of $M$ is defined to be
	\[
	M^{\prime} = \bigoplus_{\lambda \in \C} \ab\big(M_{(\lambda)})^* \subset M^*.
	\]
	Let $Y_{M^\prime}(-, x)$ denote the restriction of $ Y_{M^*}(-, x) $ to $ M^{\prime}$.
	Then $(M', Y_{M'})$ is itself a restricted $V$-module. Moreover, the canonical evaluation map yields a $V$-module isomorphism $M \cong (M^\prime)^\prime$.
	For a homomorphism $f\colon M \to N$ of restricted $V$-modules, write $f^\prime\colon N^\prime \to M^\prime$ for its restricted transpose.
	
	\subsection{Modules over tensor product VOAs}
	
	Let $V$ and $U$ be vertex operator algebras, and let $X$ be a weak $V\otimes U$-module.
	By restricting the module action to the first tensor factor $V\otimes\mathbf{1}$, $X$ inherits a weak $V$-module structure.
	Similarly, restricting the action to the second tensor factor  $\mathbf{1}\otimes U$ equips $X$ with a weak $U$-module structure.
	For later use, we set
	\[
	L^{V}(0) = (\omega_V \otimes \mathbf{1})_1, \quad L^{U}(0) = (\mathbf{1} \otimes \omega_U)_1.
	\]
	The operators $ L^V(0) $ and $ L^{U}(0)$ commute and satisfy $L(0) = L^V(0) + L^U(0).$
	
	\begin{lemma} \cite[Proposition 2.10]{McR23} \label{lem:tensor-product-hom}
		Let $V$ and $U$ be vertex operator algebras. Let $M_1, M_2$ be restricted $V$-modules and
		$N_1, N_2$ restricted $U$-modules.  Assume that $\Hom_V(M_1, M_2)$
		and $\Hom_U(N_1, N_2)$ are both finite-dimensional. Then the natural linear map
		\[
		\tau \colon \Hom_V(M_1, M_2) \otimes \Hom_U(N_1, N_2) \longrightarrow
		\Hom_{V \otimes U} \ab\big(M_1 \otimes N_1, M_2 \otimes N_2)
		\]
		defined by
		\[
		\tau(f \otimes g) (a_1 \otimes b_1) = f(a_1) \otimes g(b_1), \quad \forall \; f, g, a_1, b_1,
		\]
		is a vector space isomorphism.
	\end{lemma}
	
	\begin{lemma} \label{lem:cyclic-tensor-presentation}
		Let $V$ and $U$ be vertex operator algebras. Let $X$ be a weak $V\otimes U$-module, and let $m \in X$.
		Let $M$ (resp.\ $N$) be the cyclic $V$-submodule (resp.\ $U$-submodule) of $X$ generated by $m$.
		Then there exists a $V \otimes U$-homomorphism $f \colon M \otimes N \to X$ such that $f(m\otimes m) = m$.
	\end{lemma}
	
	\begin{proof}
		
		This lemma is obtained in the proof of Proposition 3.2 in \cite{McR23}. We present an alternative proof here.
		
		Note that $\Hom_U(N, X)$ admits a natural $V$-module structure given by
		\[
		(v_n g) (b) = v_n g(b), \qquad \forall \; g \in \Hom_U(N, X), \, v \in V, \, n \in \Z, \, b \in N.
		\]
		Let $\iota_N\colon N \hookrightarrow X$ denote the inclusion map. Then $\iota_N \in \Hom_U(N, X)$.
		The canonical evaluation map
		\[
		\theta \colon \Hom_U(N, X) \otimes N \longrightarrow X, \qquad g \otimes b \longmapsto g(b)
		\]
		is a $V\otimes U$-module homomorphism and satisfies $\theta(\iota_N \otimes m) = m$.
		On the other hand, since the $U$-module $N$ is generated by $m$, the evaluation map
		\[
		\sigma \colon \Hom_U(N, X) \longrightarrow X, \qquad g \longmapsto g(m)
		\]
		is an injective $V$-homomorphism and satisfies $\sigma(\iota_N) = m$.
		Hence there exists an injective $V$-homomorphism $r \colon M \to \Hom_U(N, X)$ such that $r(m) = \iota_N$.
		Let
		\[
		f = \theta \circ (r\otimes\id_N) \colon M \otimes N \longrightarrow X.
		\]
		Then the result follows.
	\end{proof}
	
	\begin{remark}
		
		For later arguments, it is significant that the $V$-module $M$ and the $U$-module $N$ from Lemma~\ref{lem:cyclic-tensor-presentation} are realized, respectively, as a $V$-submodule and a $U$-submodule of $X$.
		
	\end{remark}
	
	\subsection{Logarithmic intertwining operators}
	
	The definition of (logarithmic) intertwining operators and their basic properties are standard,
	and can be found in \cite{HLZ10II, M02, YZ26}.
	
	\begin{definition}\label{def:logarithmic-intertwiner}
		Let $V$ be a vertex operator algebra, and let $M, N$ and $W$ be weak $V$-modules. A \textbf{logarithmic intertwining operator} of type $\binom{W}{M \; N}$ is a linear map
		\[
		\mathcal{Y}(-, y) \colon M \longrightarrow \Hom_{\C} \ab\big(N, \; W[\log y] \{y\}), \quad
		\Y(a, y) = \sum_{n \in \C} \sum_{ k \in \N} a^{\Y}_{n; k} (\log y)^k y^{-n-1},
		\]
		where $a^{\Y}_{n; k} \in \Hom_{\C}(N, W)$,  satisfying the following conditions:
		
		\begin{enumerate}[label=(\arabic*),ref=\arabic*]
			
			\item \textbf{Truncation property.} Given any $a \in M$, $b \in N$, and $h \in \C$, there exists $s \in \N$ such that $a_{h+n;k}b = 0$ for all $n > s$ and all $k \in \N$.
			
			\item \textbf{Jacobi identity.} For all $u \in V$, $a \in M$, and $b \in N$,
			\begin{align*} \label{eq:intertwiner-jacobi}
				z^{-1} \delta\! \ab(\frac{x-y}{z}) & Y_{W}(u, x) \mathcal{Y}(a, y) b
				- z^{-1} \delta\! \ab(\frac{-y+x}{z}) \mathcal{Y}(a, y) Y_{N}(u, x) b \notag \\
				& = y^{-1} \delta\! \ab(\frac{x-z}{y}) \mathcal{Y} \ab\big(Y_{M}(u, z) a, y) b.
			\end{align*}
			
			\item \textbf{$L(-1)$-derivative property.} For all $a \in M$,
			\begin{equation*} \label{eq:intertwiner-derivative}
				\mathcal{Y} \ab\big(L(-1) a, y) = \odv{y}\mathcal{Y}(a, y).
			\end{equation*}
			
		\end{enumerate}
		
	\end{definition}
	
	For simplicity, we write $a_{n;k}b$ instead of $a^\Y_{n;k}b$ when no confusion arises.
	From now on, "intertwining operator" means "logarithmic intertwining operator" unless otherwise specified.
	Let $\mathcal{I}\binom{W}{M\;N}$ denote the space of intertwining operators of type $\binom{W}{M\;N}$.
	
	\begin{lemma}
		Let $V$ be a vertex operator algebra, and let $M$, $N$ and $W$ be weak $V$-modules, and let $\mathcal{Y}$ be an intertwining operator of type $\binom{W}{M\;N}$.
		For any $a \in M$, $b \in N$, $u \in V$, suppose that $p, q \in \N$ satisfy
		\[
		z^q Y_{M}(u, z) a \in M\dbrack{z}, \qquad x^p Y_{N}(u, x) b \in N\dbrack{x}.
		\]
		Then the following \emph{weak commutativity} and \emph{weak associativity} hold:
		\begin{enumerate}[label=(\arabic*),ref=\arabic*]
			\item \textbf{Weak commutativity}:
			\[
			(x-y)^q Y_{W}(u, x) \mathcal{Y}(a, y) = (x-y)^q \mathcal{Y}(a, y) Y_{N}(u, x).
			\]
			\item \textbf{Weak associativity}:
			\[
			(z+y)^p Y_{W}(u, z+y) \mathcal{Y}(a, y) b = (z+y)^p \mathcal{Y} \ab\big(Y_{M}(u, z) a, \, y) b.
			\]
		\end{enumerate}
	\end{lemma}
	
	\vspace{1\baselineskip}
	
	\noindent \textbf{ The operators $\Omega_r$ and $T_r$}.
	We recall the definitions and properties of $\Omega_r\Y$ and $T_r\Y$ from \cite{HLZ10II, FHL93}.
	We point out that in this paper we use the notation $T_r\Y$ instead of $A_r\Y$ in \cite{HLZ10II}.
	Let $r \in \Z$. Let $M, N$ and $W$ be weak $V$-modules over a vertex operator algebra $V$, and
	let $\Y$ be an intertwining operator of type $\binom{W}{M \; N}$. Define a linear map
	\[
	\Omega_r \Y \colon N \longrightarrow \Hom_{\C} \ab\big(M, \; W [\log y] \{y\})
	\]
	by
	\begin{equation} \label{eq:skew-symmetry-transform}
		(\Omega_r \Y) \ab\big(b, y) a = e^{yL(-1)} \Y \ab\big(a, e^{ (2r+1) \pi \mathrm{i}}y) b,
	\end{equation}
	for all $ a \in M, b \in N.$
	Then $\Omega_r\Y$ is an intertwining operator of type  $\binom{W}{N \; M}$.
	Moreover, it satisfies
	\[
	\Omega_{-r-1} (\Omega_r \Y) = \Omega_{r} (\Omega_{-r-1}\Y) = \Y.
	\]
	
	Furthermore, assume that $M$, $N$, and $W$ are restricted $V$-modules. Define a linear map
	\[
	T_r\Y\colon M \longrightarrow \Hom_{\C} \ab\big(W^\prime, \; N^\prime [\log y] \{y\})
	\]
	by
	\begin{equation}\label{eq:contragredient-transform}
		\ab\big\langle(T_r\Y) (a, y) w^{\prime}, \, b\rangle
		= \ab\big\langle w^{\prime}, \; \Y \ab\big(e^{yL(1)} \ab\big(e^{(2r+1) \pi\mathrm{i}} y^{-2})^{L(0)} a, \, y^{-1}) b \rangle
	\end{equation}
	for all $a \in M$, $b \in N$, and $w' \in W'$. Then $T_r\Y$ is an intertwining operator of type $\binom{N'}{M, \; W'}$,
	and satisfies
	\[
	T_{-r-1} \ab\big(T_r\Y) = T_{r} \ab\big(T_{-r-1}\Y) = \Y.
	\]
	
	\vspace{1\baselineskip}
	
	\noindent \textbf{ The trifunctor $\mathcal{I} \binom{-}{- \; -}$}. Let $M, M_0, N, N_0$ and $W, W_0$ be weak $V$-modules over a vertex operator algebra $V$,
	and let
	\[
	h \colon M_0 \longrightarrow M, \; g \colon N_0 \longrightarrow N, \; f \colon W \longrightarrow W_0
	\]
	be $V$-module homomorphisms.
	Let $\mathcal{Y}(-, x)$ be an intertwining operator of type $\binom{W}{M \; N}$.
	We define a new operator
	\[
	\ab\big([h, g, f] \circ \mathcal{Y}) (-, y) \colon M_0 \longrightarrow \Hom_{\C} \ab\big(N_0, \; W_0 [\log y] \{y\})
	\]
	by the formula
	\begin{equation}
		\ab\big([h, g, f] \circ \mathcal{Y}) (a_0, y) b_0 := f \ab\big(\mathcal{Y}(ha_0, y) g b_0)
	\end{equation}
	for all $a_0 \in M_0$ and $b_0 \in N_0$.
	By \cite[Lemma 3.6]{YZ26}, the following statements hold.
	First, $ [h, g, f] \circ \mathcal{Y} $  is an intertwining operator of type $ \binom{W_0}{M_0 \; N_0}$.
	For any $r \in \Z$, the following identity holds:
	\begin{equation}
		\Omega_r \ab\big([h, g, f] \circ \mathcal{Y})
		=
		[g, h, f] \circ(\Omega_r\mathcal{Y}).
		\label{eq:skew-symmetry-naturality}
	\end{equation}
	If we further assume that these modules are restricted $V$-modules, then we have
	\begin{equation}
		T_r \ab\big([h, g, f] \circ\mathcal{Y}) = [h, f^\prime, g^\prime] \circ (T_r \mathcal{Y}).
		\label{eq:contragredient-naturality}
	\end{equation}
	
	Now consider the product category
	\[
	(\mathfrak{M}^w_V)^{\opp} \times (\mathfrak{M}^w_V)^{\opp} \times \mathfrak{M}^w_V.
	\]
	Its objects are triples $(M, N, W)$ of objects in $\mathfrak{M}^w_V$.  A morphism from $(M, N, W)$ to $(M_0, N_0, W_0)$ is a triple $(h, g, f)$ of $V$-module homomorphisms
	$h\colon M_0 \to M$, $g \colon N_0 \to N$, $f \colon W \to W_0$, with morphisms composed componentwise.
	We define the trifunctor
	\[
	\mathcal{I}: (\mathfrak{M}^w_V)^{\opp} \times (\mathfrak{M}^w_V)^{\opp} \times \mathfrak{M}^w_V \longrightarrow \mathbf{Vect},
	\]
	as follows:
	\begin{enumerate}[label=(\arabic*),ref=\arabic*]
		\item On objects: $\mathcal{I}(M, N, W)$ is the vector space of intertwining operators of type $\mathcal{I}\binom{W}{M \; N}$.
		
		\item On morphisms: each morphism $(h, g, f)$ induces a linear map
		\[
		\mathcal{I}(h, g, f) \colon \mathcal{I}\binom{W}{M, N} \longrightarrow \mathcal{I}\binom{W_0}{M_0, N_0}, \quad \Y \longrightarrow [h, g, f] \circ \Y.
		\]
	\end{enumerate}
	Functoriality of $\mathcal{I}$ holds trivially from the componentwise composition of morphisms in the product category.
	Given a fixed weak $V$-module $W$, we obtain a bifunctor
	\[
	\mathcal{I}^{W}(-, -) := \mathcal{I}(-, -, W) \colon (\mathfrak{M}_V^w)^{\opp} \times (\mathfrak{M}_V^w)^{\opp} \to \mathbf{Vect}.
	\]
	Note that the bifunctor $\mathcal{I}^{W}$ is bilinear and left exact in each variable.

	\section{Constructing \texorpdfstring{$\mathcal{R}(W)$}{R(W)} via Intertwining Operators} \label{sec:regular-construction}
	
	Let $V$ be a vertex operator algebra. Let $W$ be a weak $V$-module.
	Define
	\begin{equation} \label{eq:regular-construction}
		\mathcal{R}(W)
		=
		\Span_{\C}
		\ab\{
		\mathcal{Y}(a, y) b
		\; \middle| \;
		M, N \in \mathfrak{M}^w_V, \;
		a \in M, \; b \in N, \;
		\mathcal{Y} \in \mathcal{I}\binom{W}{M\; N}
		\}
		\subset W[\log y] \{y\}.
	\end{equation}
	
	Consider an element
	\[
	\sum_{i = 1}^n \lambda_i \mathcal{Y}_i(a_i, y) b_i \in \mathcal{R}(W),
	\]
	where for each $i$, $\lambda_i \in \C$, $M_i, N_i \in \mathfrak{M}^w_V$, $a_i \in M_i$, $b_i \in N_i$, and
	$\mathcal{Y}_i \in \mathcal{I}\binom{W}{M_i\;N_i}$.
	Set
	\begin{align*}
		M & = M_1 \oplus \cdots \oplus M_n, \quad
		N = N_1 \oplus \cdots \oplus N_n, \quad
		\mathcal{Y} = \mathcal{Y}_1 \oplus \cdots \oplus \mathcal{Y}_n, \\
		a & = (\lambda_1 a_1, \ldots, \lambda_n a_n) \in M, \quad \text{and} \quad
		b = (b_1, \ldots, b_n) \in N.
	\end{align*}
	Here the direct-sum notation is defined by
	\[
	\mathcal{Y} \ab\big((a_i)_i, y) (b_i)_i
	:= \sum_{i=1}^n \mathcal{Y}_i(a_i, y) b_i.
	\]
	Then $\mathcal{Y}$ is an intertwining operator of type $\binom{W}{M\;N}$ such that
	\[
	\sum_{i = 1}^n \lambda_i \mathcal{Y}_i(a_i, y) b_i = \mathcal{Y}(a, y) b \in \mathcal{R}(W).
	\]
	It follows that every element of $\mathcal{R}(W)$ can be expressed as $\mathcal{Y}(a, y) b$ for some weak $V$-modules $M$ and $N$, elements $a \in M$ and $b \in N$, and an intertwining operator $\mathcal{Y} \in \mathcal{I}\binom{W}{M\;N}$.
	For notational simplicity, in what follows, when we write
	\[
	\mathcal{Y}(a, y) b \in \mathcal{R}(W),
	\]
	we implicitly mean that $\mathcal{Y}$ is an intertwining operator of type
	$\binom{W}{M\,N}$ for some weak $V$-modules $M$ and $N$, with
	$a \in M$ and $b \in N$.
	
	For $u \in V$, we define two linear maps
	\[
	Y_{\mathcal{R}(W)}^L(u, x), \quad Y_{\mathcal{R}(W)}^R(u, x)
	\colon \mathcal{R}(W) \longrightarrow \mathcal{R}(W) \dparen{x}
	\]
	by the formulas
	\begin{align} \label{eq:regular-left-right-actions}
		Y_{\mathcal{R}(W)}^L(u, x) \,\mathcal{Y}(a, y) b
		& = \mathcal{Y} \ab\big(Y_{M}(u, x) a, \,y) b, \\
		Y_{\mathcal{R}(W)}^R(u, x) \,\mathcal{Y}(a, y) b
		& = \mathcal{Y}(a, y) Y_{N}(u, x) b,
	\end{align}
	for all $\mathcal{Y}(a, y) b \in \mathcal{R}(W)$, and extend linearly.
	
	\begin{proposition}
		Let $V$ be a vertex operator algebra and $W$ a weak $V$-module.
		For all $u, v \in V$, the operators
		\[
		Y_{\mathcal{R}(W)}^L(v, x), \quad Y_{\mathcal{R}(W)}^R(u, x)
		\colon \mathcal{R}(W) \longrightarrow \mathcal{R}(W) \dparen{x}
		\]
		are well-defined and mutually commute.
		Moreover, both $(\mathcal{R}(W), Y_{\mathcal{R}(W)}^L)$ and $(\mathcal{R}(W), Y_{\mathcal{R}(W)}^R)$ are weak $V$-modules.
	\end{proposition}
	
	\begin{proof}
		
		Let $u, v \in V$. We first prove that the maps $Y_{\mathcal{R}(W)}^R(u, x)$ and $Y_{\mathcal{R}(W)}^L(v, x)$ are both well-defined.
		Suppose that
		\[
		\sum_{i = 1}^n \lambda_i \mathcal{Y}_i(a_i, y) b_i = 0,
		\]
		where for each $i$, $\lambda_i \in \C$, and $\mathcal{Y}_i$ is an intertwining operator of type $\binom{W}{M_i \; N_i}$ for some weak $V$-modules $M_i$ and $N_i$, with $a_i \in M_i$ and $b_i \in N_i$.
		
		By the weak commutativity of $\mathcal{Y}_i$, there exists $k_i \in \N$ such that
		\[
		(x-y)^{k_i} \mathcal{Y}_i(a_i, y) Y_{N_i}(u, x) b_i
		= (x-y)^{k_i} Y_{W}(u, x) \mathcal{Y}_i(a_i, y) b_i.
		\]
		Let $k \in \N$ satisfy $k > \max \{k_1, k_2, \dots, k_n\}$. Then
		\[
		\sum_{i = 1}^n (x-y)^k \lambda_i \mathcal{Y}_i(a_i, y) Y_{N_i}(u, x) b_i
		= \sum_{i = 1}^n (x-y)^k \lambda_i Y_{W}(u, x) \mathcal{Y}_i(a_i, y) b_i
		= 0.
		\]
		Since the first expression vanishes, multiplying it by $(-y+x)^{-k}$ yields
		\[
		\sum_{i = 1}^n \lambda_i \mathcal{Y}_i(a_i, y) Y_{N_i}(u, x) b_i = 0.
		\]
		From the definition of $Y_{\mathcal{R}(W)}^R$, we therefore obtain
		\[
		Y^R_{\mathcal{R}(W)}(u, x) \ab(\sum_{i = 1}^n \lambda_i \mathcal{Y}_i(a_i, y) b_i)
		= \sum_{i = 1}^n \lambda_i \mathcal{Y}_i(a_i, y) Y_{N_i}(u, x) b_i = 0.
		\]
		This proves that $Y_{\mathcal{R}(W)}^R$ is well-defined.
		
		Similarly, by the weak associativity of $\mathcal{Y}_i$, there exists a sufficiently large integer $k \in \N$ such that
		\[
		\sum_{i = 1}^n (z+y)^k \lambda_i \mathcal{Y}_i \ab\big(Y_{M_i}(u, z) a_i, y) b_i
		= \sum_{i = 1}^n (z+y)^k \lambda_i Y_{W}(u, z+y) \mathcal{Y}_i(a_i, y) b_i
		= 0.
		\]
		Multiplying the first expression by $(y+z)^{-k}$ and using the definition of $Y_{\mathcal{R}(W)}^L$, we get
		\[
		Y^L_{\mathcal{R}(W)}(u, x) \ab(\sum_{i = 1}^n \lambda_i \mathcal{Y}_i(a_i, y) b_i)
		= \sum_{i = 1}^n \lambda_i \mathcal{Y}_i \ab\big(Y_{M_i}(u, x) a_i, y) b_i
		= 0.
		\]
		Hence the operator $Y_{\mathcal{R}(W)}^L$ is also well-defined.
		
		From the definitions of $Y_{\mathcal{R}(W)}^R$ and $Y_{\mathcal{R}(W)}^L$,
		it is immediate that $Y_{\mathcal{R}(W)}^R(u, x)$ and $Y_{\mathcal{R}(W)}^L(v, x)$ commute for all $u, v \in V$.
		
		Next, we prove that the pair $\ab\big(\mathcal{R}(W), Y^L_{\mathcal{R}(W)})$ forms a weak $V$-module.
		Let $\mathcal{Y}(a, y) b \in \mathcal{R}(W)$, where $\mathcal{Y}$ is an intertwining operator of type $\binom{W}{M\;N}$
		for some weak $V$-modules $M, N$ with $a \in M$ and $b \in N$.
		By the definition of $Y^L_{\mathcal{R}(W)}$ and the Jacobi identity for $Y_{M}$, for all $u, v \in V$ we have
		\begin{align*}
			& x_0^{-1}\delta \ab(\frac{x_1-x_2}{x_0}) Y_{\mathcal{R}(W)}^L(u, x_1) Y_{\mathcal{R}(W)}^L(v, x_2) \mathcal{Y}(a, y) b \\
			& \quad - x_0^{-1}\delta \ab(\frac{-x_2+x_1}{x_0}) Y_{\mathcal{R}(W)}^L(v, x_2) Y_{\mathcal{R}(W)}^L(u, x_1) \mathcal{Y}(a, y) b \\
			& = x_0^{-1}\delta \ab(\frac{x_1-x_2}{x_0}) \mathcal{Y} \ab\bigg(Y_{M}(u, x_1) Y_{M}(v, x_2) a, \, y) b \\
			& \quad - x_0^{-1}\delta \ab(\frac{-x_2+x_1}{x_0}) \mathcal{Y} \ab\bigg(Y_{M}(v, x_2) Y_{M}(u, x_1) a, \, y) b \\
			& = x_2^{-1}\delta \ab(\frac{x_1-x_0}{x_2}) \mathcal{Y} \ab\bigg(Y_{M} \ab\big(Y(u, x_0) v, x_2) a, \, y) b \\
			& = x_2^{-1}\delta \ab(\frac{x_1-x_0}{x_2}) Y_{\mathcal{R}(W)}^L \ab\big(Y(u, x_0) v, x_2) \mathcal{Y}(a, y) b.
		\end{align*}
		Thus the Jacobi identity holds for $Y^L_{\mathcal{R}(W)}$.
		The vacuum property and the lower-truncation property are straightforward to verify.
		Consequently, $\ab\big(\mathcal{R}(W), Y^L_{\mathcal{R}(W)})$ is a weak $V$-module.
		
		Similarly, using the definition of $Y_{\mathcal{R}(W)}^R$ and the Jacobi identity for $Y_{N}$,
		we conclude that $\ab\big(\mathcal{R}(W), Y_{\mathcal{R}(W)}^R)$ also forms a weak $V$-module.
		This completes the proof.
		
	\end{proof}
	
	Let $V$ be a vertex operator algebra, and let $W$ be a weak $V$-module.
	Now we can define a well-defined linear map
	\[
	Y_{\mathcal{R}(W)}(-, x) -: (V \otimes V) \otimes \mathcal{R}(W) \to \mathcal{R}(W) \dparen{x}
	\]
	by
	\[
	Y_{\mathcal{R}(W)}(u \otimes v, x) = Y_{\mathcal{R}(W)}^L(u, x) Y_{ \mathcal{R}(W)}^R(v, x), \quad \forall \; u, v \in V.
	\]
	Since $Y_{\mathcal{R}(W)}^L(u, x)$ and $Y_{\mathcal{R}(W)}^R(v, x)$ action on $\mathcal{R}(W)$ commute for any $u, v \in V$, we immediately have
	
	\begin{theorem}
		Let $V$ be a vertex operator algebra, and let $W$ be a weak $V$-module.
		Then the pair $\ab\big(\mathcal{R}(W), Y_{\mathcal{R}(W)})$ carries the structure of a weak $V \otimes V$-module.
	\end{theorem}
	
	Let $V$ be a vertex operator algebra, and let $f \colon W_1 \to W_2$ be a homomorphism of weak $V$-modules.
	Then the map
	\[
	\mathcal{R}(f) \colon \mathcal{R}(W_1) \longrightarrow \mathcal{R}(W_2)
	\]
	defined by
	\[
	\mathcal{R}(f) \ab\big(\mathcal{Y}(a, y) b) = f \ab\big(\mathcal{Y}(a, y) b)
	= \ab\big([\id_{M}, \id_{N}, f] \circ \mathcal{Y}) (a, y) b \in \mathcal{R}(W_2)
	\]
	for all $\Y(a, y) b \in \mathcal{R}(W_1)$,
	is well-defined and is a $V \otimes V$-module homomorphism.
	Thus $\mathcal{R}$ becomes a functor from the category of weak $V$-modules to the category of weak $V\otimes V$-modules.
	The following proposition show that the functor $\mathcal{R}$ is left-exact.
	
	\begin{proposition}
		Let $ W_1, W_2, W_3$ be weak $V$-modules over a vertex operator algebra $V$, and suppose
		\[
		0 \to W_1 \xrightarrow{f} W_2 \xrightarrow{g} W_3 \to 0
		\]
		is an exact sequence of $V$-module homomorphisms.
		Then, the induced sequence of weak $V \otimes V$-modules
		\[
		0 \to \mathcal{R}(W_1) \xrightarrow{\mathcal{R}(f)} \mathcal{R}(W_2) \xrightarrow{\mathcal{R}(g)} \mathcal{R}(W_3)
		\]
		is exact.
		
	\end{proposition}
	
	\begin{proof}
		It is straightforward to verify that $\mathcal{R}(f)$ is injective and that
		\[
		\Image\mathcal{R}(f) \subseteq \Ker\mathcal{R}(g).
		\]
		
		Let $\mathcal{Y}(a, y) b \in \Ker\mathcal{R}(g)$,
		where $\mathcal{Y}$ is an intertwining operator of type $\binom{W_2}{M\,N}$ for some weak $V$-modules $M$ and $N$, with $a \in M$ and $b \in N$.
		Without loss of generality, we may assume that $M$ is generated by $a$ and $N$ is generated by $b$; indeed, one may replace $M$ and $N$ by the submodules generated by $a$ and $b$, respectively.
		Then $g \ab\big(\mathcal{Y}(a, y) b) = 0$.
		It follows that
		\begin{align*}
			g \ab\big(\mathcal{Y}(u_n a, y) v_m b)
			& = \ab\big([\id_{M}, \id_{N}, g] \circ \mathcal{Y}) (u_n a, y) v_m b \\
			& = u^L_n v^R_m \ab\big(\ab\big([\id_{M}, \id_{N}, g] \circ \mathcal{Y}) (a, y) b) \\
			& = u^L_n v^R_m \ab\big(g \ab\big(\mathcal{Y}(a, y) b)) = 0.
		\end{align*}
		Consequently,
		\[
		\Image\mathcal{Y}
		= \Span_{\C} \ab\{ c_{n;k}d \; \middle| \; c \in M, \ d \in N, \ n \in \C, \ k \in \N \}
		\subseteq \Ker g = \Image f.
		\]
		Let $f^{-1}\colon \Image f \to W_1$ denote the inverse of the induced linear isomorphism $f\colon W_1 \to \Image f$.
		Then $[\id_{M}, \id_{N}, f^{-1}] \circ \mathcal{Y}$ is an intertwining operator of type $\binom{W_1}{M\,N}$ and satisfies
		\[
		\mathcal{R}(f) \ab\big(\ab\big([\id_{M}, \id_{N}, f^{-1}] \circ \mathcal{Y}) (a, y) b) = \mathcal{Y}(a, y) b.
		\]
		It follows that $\Ker\mathcal{R}(g) \subseteq \Image\mathcal{R}(f)$.
		This completes the proof.
	\end{proof}
	
	\begin{lemma}
		Let $V$ be a vertex operator algebra, and let $W$ be a weak $V$-module.
		Let $\phi(y) \in \mathcal{R}(W)$. Then
		\begin{align} \label{eq:regular-translation}
			\omega^L_0 \phi(y) = \odv{y} \phi(y).
		\end{align}
		Moreover the following identity holds for all $u \in V$:
		\begin{align} \label{eq:mixed-jacobi}
			& z^{-1} \delta \ab(\frac{x-y}{z}) Y_{W}(u, x) \phi(y)
			-z^{-1} \delta \ab(\frac{-y+x}{z}) Y^R_{\mathcal{R}(W)}(u, x) \phi(y) \notag \\
			& = y^{-1} \delta \ab(\frac{x-z}{y}) Y^L_{\mathcal{R}(W)}(u, z) \phi(y).
		\end{align}
	\end{lemma}
	
	\begin{proof}
		Since $\phi(y) = \Y(a, y) b$ for some intertwining operator $ \Y \in \mathcal{I}\binom{W}{M \; N}$,
		the definitions of $Y^L_{\mathcal{R}(W)}$ and $ Y^R_{\mathcal{R}(W)}$
		imply that identities \eqref{eq:regular-translation} and \eqref{eq:mixed-jacobi} are respectively equivalent to the $L(-1)$-derivative property and the Jacobi identity of $\Y$.
		
	\end{proof}
	
	\begin{theorem} \label{thm:hom-realization}
		Let $V$ be a vertex operatro algebra, and let $W$ be a weak $V$-module. For any weak $V$-modules $M$ and $N$, the linear map
		\begin{align}\label{eq:hom-realization-map}
			\psi^{W}_{M, N}\colon \mathcal{I} \binom{W}{M \; N} \longrightarrow \Hom_{V \otimes V} \ab\big(M \otimes N, \; \mathcal{R}(W))
		\end{align}
		defined by
		\[
		\psi^{W}_{M, N}(\Y) (a \otimes b) = \Y(a, y) b \quad \forall \; a \in M, b \in N,
		\]
		is a linear isomorphism.  Moreover, these isomorphisms assemble into a natural isomorphism
		\[
		\psi^{W}\colon \mathcal{I}^{W}(-, -) \Rightarrow  
		\Hom_{V \otimes V} \ab\big(- \otimes -, \; \mathcal{R}(W)) := \mathcal{H}_{\mathcal{R}(W)},
		\]
		where both are bifunctors from $(\mathfrak{M}^w_V)^{\opp}\times(\mathfrak{M}^w_V)^{\opp}$ to $\mathbf{Vect}$.
		
	\end{theorem}
	
	\begin{proof}
		
		Let $M$ and $N$ be weak $V$-modules. The injectivity of $\psi^W_{M, N}$ is immediate.
		For each $\mathcal{Y} \in \mathcal{I}\binom{W}{M\;N}$, the map
		$\psi^W_{M, N}(\mathcal{Y}) \colon M\otimes N \to \mathcal{R}(W)$ is a
		$V\otimes V$-module homomorphism. To establish surjectivity, take an arbitrary
		\[
		\mathcal{J} \in \Hom_{V \otimes V} \ab\big(M \otimes N, \; \mathcal{R}(W)).
		\]
		By the canonical isomorphism
		\[
		\Hom_{\C} \ab\big(M \otimes N, \; W) \cong \Hom_{\C} \ab\big(M, \; \Hom_{\C}(N, W)),
		\]
		For any $a \in M$ and $b \in N$,  we may express $\mathcal{J}(a\otimes b) $ in the form
		\[
		\mathcal{J}(a, y) b = \sum_{n \in \C} \sum_{k = 0}^{p_{a, b, n}} a_{n;k}b\, (\log y)^k y^{-n-1},
		\]
		where $p_{a, b, n} \in \N$ and $a_{n;k} \in \Hom_{\C}(N, W)$.
		We will show that $\mathcal{J}(-, y)$ defines an intertwining operator of type $\binom{W}{M\;N}$ and that $\psi^W_{M, N} \ab\big(\mathcal{J}(-, y)) = \mathcal{J}$.
		
		Since $\mathcal{J}(a, y) b \in \mathcal{R}(W)$ for all $a, b$, the truncation property holds for $\mathcal{J}(-, y)$.
		Because $\mathcal{J}$ is a $V\otimes V$-module homomorphism, we apply identity \eqref{eq:regular-translation} to verify the $L(-1)$‑derivative property for $\mathcal{J}(-, y)$:
		\begin{align*}
			\mathcal{J} \ab\big(L(-1) a, y) b
			& = \mathcal{J} \ab\big(L(-1) a \otimes b) \\
			& = \omega^L_0 \mathcal{J}(a\otimes b)
			= \odv{y}\mathcal{J}(a\otimes b)
			= \odv{y}\mathcal{J}(a, y) b,
			\quad \forall\, a \in M, \, b \in N.
		\end{align*}
		Furthermore, applying identity \eqref{eq:mixed-jacobi} and the fact that $\mathcal{J}$ is a $V\otimes V$-module homomorphism, we obtain:
		\begin{align*}
			& z^{-1} \delta \ab(\frac{x-y}{z}) Y_{W}(u, x) \mathcal{J}(a, y) b
			- z^{-1} \delta \ab(\frac{-y+x}{z}) \mathcal{J}(a, y) Y_{N}(u, x) b \\
			& = z^{-1} \delta \ab(\frac{x-y}{z}) Y_{W}(u, x) \mathcal{J}(a\otimes b)
			- z^{-1} \delta \ab(\frac{-y+x}{z}) \mathcal{J} \ab\big(a\otimes Y_{N}(u, x) b) \\
			& = z^{-1} \delta \ab(\frac{x-y}{z}) Y_{W}(u, x) \mathcal{J}(a\otimes b)
			- z^{-1} \delta \ab(\frac{-y+x}{z}) Y^{R}_{\mathcal{R}(W)}(u, x) \mathcal{J}(a\otimes b) \\
			& = y^{-1} \delta \ab(\frac{x-z}{y}) Y^L_{\mathcal{R}(W)}(u, z) \mathcal{J}(a\otimes b) \\
			& = y^{-1} \delta \ab(\frac{x-z}{y}) \mathcal{J} \ab\big(Y_{M}(u, z) a \otimes b) \\
			& = y^{-1} \delta \ab(\frac{x-z}{y}) \mathcal{J} \ab\big(Y_{M}(u, z) a, y) b,
			\quad \forall\, u \in V, \, a \in M, \, b \in N.
		\end{align*}
		This proves that $\mathcal{J}(-, y)$ satisfies the Jacobi identity.
		Therefore, $\mathcal{J}(-, y)$ is an intertwining operator of type $\binom{W}{M \; N}$, and $\psi^W_{M, N} \ab\big(\mathcal{J}(-, y)) = \mathcal{J}$.
		Hence $\psi^W_{M, N}$ is surjective.
		
		\vspace{0.3cm}
		
		The functoriality of $\psi^W$ amounts to the assertion that for any weak $V$-module homomorphisms $f\colon M_0 \to M$ and $g\colon N_0 \to N$, the diagram below commutes:
		\[
		\begin{tikzcd}
			\mathcal{I}\binom{W}{M,N}
			\arrow[r,"{\psi^W_{M,N}}"]
			\arrow[d,"{[f,g,\id_{W}] \circ -}"']
			& \Hom_{V\otimes V} \ab\big(M\otimes N,\mathcal{R}(W))
			\arrow[d,"{\Hom_{V\otimes V}(f\otimes g,\mathcal{R}(W))}"] \\
			\mathcal{I}\binom{W}{M_0,N_0}
			\arrow[r,"{\psi^W_{M_0,N_0}}"']
			& \Hom_{V\otimes V} \ab\big(M_0\otimes N_0,\mathcal{R}(W))
		\end{tikzcd}
		\]
		We verify commutativity by the following chain of equalities, which holds for all $a_0 \in M_0$, $b_0 \in N_0$, and $\mathcal{Y} \in \mathcal{I}\binom{W}{M\;N}$:
		\begin{align*}
			& \ab\big(\Hom_{V\otimes V} \ab\big(f\otimes g, \mathcal{R}(W)) \circ \psi^W_{M, N}(\mathcal{Y})) (a_0\otimes b_0) \\
			& \quad = \psi^W_{M, N}(\mathcal{Y}) \ab\big((f\otimes g) (a_0\otimes b_0)) \\
			& \quad = \mathcal{Y}(fa_0, y) gb_0 = \ab\big([f, g, \id_{W}] \circ \mathcal{Y}) (a_0, y) b_0 \\
			& \quad = \psi^W_{M_0, N_0} \ab\big([f, g, \id_{W}] \circ \mathcal{Y}) (a_0\otimes b_0).
		\end{align*}
		This completes the proof.
		
	\end{proof}
	
	\section{ The \texorpdfstring{$C_2$}{C₂}-Cofinite Case: Finiteness and The Deligne Tensor Product Approach} \label{sec:c2-cofinite-case}
	
	
	Let $V$ be a $C_2$-cofinite vertex operator algebra with nonnegative integer weights.
	The category of restricted $V$-modules enjoys the following properties (see \cite{ABD04, DLM00, M04, H09, YZ26}):
	\begin{enumerate}[label=(\arabic*),ref=\arabic*]
		\item there exist only finitely many isomorphism classes of simple restricted $V$-modules;
		\item every restricted $V$-module is of finite length, which in particular forces all $\Hom$-spaces between restricted $V$-modules to be automatically finite-dimensional over $\C$;
		\item every restricted $V$-module admits a projective cover, or equivalently, every simple restricted $V$-module admits a projective cover;
		\item the space of intertwining operators among any triple of restricted $V$-modules is finite-dimensional.
	\end{enumerate}
	Properties (1), (2), and (3) imply that the category of restricted $V$-modules is a finite abelian category.
	
	Let $S_1, S_2, \dots, S_m$ be a complete set of representatives of isomorphism classes of irreducible restricted $V$-modules.
	For each $i$, let $P_i$ be the projective cover of $S_i$, and set $P = \bigoplus_{i = 1}^m P_i$.
	Then $P$ is a projective generator of the category of restricted $V$-modules.
	Moreover, since every cyclic submodule of a weak $V$-module is a restricted $V$-module \cite{ABD04,M04},
	it follows that $P$ is also a projective generator in the category of weak $V$-modules.
	
	For convenience, we introduce the following terminology.
	Let $R$ be a ring or a vertex operator algebra, and let $\mathcal{U}$ be a class of $R$-modules. We say that $\mathcal{U}$ \textbf{generates} an $R$-module $M$, or that $M$ is \textbf{generated by $\mathcal{U}$},
	if there exists a family $(U_i)_{i \in I}$ of modules in $\mathcal{U}$ together with a surjective homomorphism
	\[
	\bigoplus_{i \in I} U_i \twoheadrightarrow M \to 0.
	\]
	If $\mathcal{U} = \{U\}$ is a singleton class, we simply say that $U$ generates $M$. This means that there exists some index set $I$ and an epimorphism
	\[
	U^{\oplus I} \twoheadrightarrow M \to 0.
	\]
	
	\begin{proposition} \label{prop:c2-regular-finiteness}
		Let $V$ be a $C_2$-cofinite vertex operator algebra with nonnegative integer weights, and let $W$ be a restricted $V$-module. Then $\mathcal{R}(W)$ is also a restricted $V \otimes V$-module. In particular, $\mathcal{R}(W)$ has finite composition length.
	\end{proposition}
	
	\begin{proof}
		Let
		\[
		\mathcal{U} = \ab\{ M \otimes N \; \big| \; M, N \in \mathfrak{M}^w_V  \}.
		\]
		It follows from Lemma~\ref{lem:cyclic-tensor-presentation} that the class $\mathcal{U}$  generates all weak $V \otimes V$-modules.
		In particular,  $\mathcal{U}$ generates $\mathcal{R}(W)$.
		Let $P$ be a projective generator of the category of restricted $V$-modules.
		Then $P$ is finitely generated, and consequently $P\otimes P$ is a restricted $V \otimes V$-module.
		Moreover, since $P$ generates every weak $V$-module, $P \otimes P$ generates every object of $\mathcal{U}$, and therefore also generates $\mathcal{R}(W)$.
		
		On the other hand, combining the finite-dimensionality of $\mathcal{I} \binom{W}{P \; P}$
		with the isomorphism \eqref{eq:hom-realization-map}, we obtain the finite-dimensionality of $\Hom_{V \otimes V} \ab\big(P \otimes P, \; \mathcal{R}(W))$.
		Choosing a basis $f_1, \cdots, f_s$ of this $\Hom$-space, we obtain
		\[
		\mathcal{R}(W) = \sum_{i = 1}^s f_i(P \otimes P).
		\]
		Thus $\mathcal{R}(W)$ is a restricted $V\otimes V$-module. This completes the proof.
		
	\end{proof}
	
	\begin{remark}
		Let $V$ be a $C_2$-cofinite vertex operator algebra with nonnegative integer weights.
		By \cite[Lemma~3.29]{GZ26}, a generalized $V$-module $W$ has finite
		composition length if and only if $\Hom_V(P, W)$ is
		finite-dimensional for every projective $V$-module $P$ of finite length.
		This condition is easily seen to be equivalent to the finite-dimensionality of
		$\Hom_V(Q, W)$ for a projective generator $Q$ of finite length.
		This observation yields an alternative proof of Proposition~\ref{prop:c2-regular-finiteness}.
		
	\end{remark}
	
	\begin{example}
		Let $V$ be a $C_2$-cofinite, rational vertex operator algebra with nonnegative integer weights.
		Let $S_1, S_2, \dots, S_m$ be a complete set of representatives of isomorphism classes of irreducible ordinary $V$-modules.
		Let $W$ be an ordinary $V$-module. Then
		\[
		\mathcal{R}(W) \cong \bigoplus_{1 \leq i, \; j \leq m} (S_i \otimes S_j)^{\oplus N(S_i, S_j; W)},
		\]
		is also a  ordinary $V\otimes V$-module,  where $N(S_i, S_j;W)$ denotes the dimension of the intertwining operator space $\mathcal{I}\binom{W}{S_i\;S_j}$.
	\end{example}
	
	\begin{proof}
		Since $V \otimes V$ is rational, Proposition~\ref{prop:c2-regular-finiteness} implies that $\mathcal{R}(W)$ is fact a ordinary $V\otimes V$-module.
		Moreover, $\mathcal{R}(W)$ is completely determined up to isomorphism by the multiplicities of the simple modules $S_i \otimes S_j$ occurring in it.
		The desired result then follows from Theorem~\ref{thm:hom-realization}.
	\end{proof}
	
	In the remainder of this section, we use the Deligne tensor product realization for finite $\C$-linear abelian categories to establish a Yoneda-type lemma for spaces of intertwining operators in the $C_2$-cofinite case. The next section extends this to the general case.
	
	Let $\mathcal{C}$ and $\mathcal{D}$ be categories, and let $F, G \colon \mathcal{C} \to \mathcal{D} $ be functors.
	Define
	\[
	\Nat(F, G) = \ab\{\text{natural transformations from $F$ to $G$}\}.
	\]
	For $\C$-linear abelian categories $\mathcal{A}, \mathcal{B}, \mathcal{C}$, let $\Rex \ab\big(\mathcal{A}, \mathcal{B})$ denote the category of $\C$-linear right exact functors from
	$\mathcal{A}$ to $\mathcal{B}$, and let $\Brex \ab\big(\mathcal{A}\times\mathcal{B}, \mathcal{C})$
	denote the category of $\C$-bilinear bifunctors from $\mathcal{A}\times\mathcal{B}$ to $\mathcal{C}$ which are right exact in each variable.

	\begin{definition} \cite{Del90}
		Let $\mathcal{A}$ and $\mathcal{B}$ be $\C$-linear abelian categories.
		A \textbf{Deligne tensor product} of $\mathcal{A}$ and $\mathcal{B}$ is a $\C$-linear abelian category $\mathcal{A}\boxtimes\mathcal{B}$, together with a $\C$-bilinear bifunctor
		\[
		\boxtimes \colon \mathcal{A}\times\mathcal{B} \longrightarrow \mathcal{A}\boxtimes\mathcal{B}
		\]
		which is right exact in each variable, such that for every $\C$-linear abelian category $\mathcal{C}$, the functor
		\[
		\Rex(\mathcal{A}\boxtimes\mathcal{B}, \; \mathcal{C})
		\longrightarrow
		\Brex(\mathcal{A}\times\mathcal{B}, \; \mathcal{C}),
		\qquad
		F \longmapsto F\circ \boxtimes,
		\]
		is an equivalence of categories.
		
	\end{definition}
	
	Note that if $\mathcal{A}$ and $\mathcal{B}$ are locally finite $\C$-linear abelian categories,
	then their Deligne tensor product $\mathcal{A} \boxtimes \mathcal{B}$ exists and is unique up to isomorphism.
	Assume that the Deligne tensor product $\mathcal{A}\boxtimes\mathcal{B}$ of $\mathcal{A}$ and $\mathcal{B}$ exists.  Then the equivalence of functor categories can be restated as follows:
	\begin{enumerate}[label=(\arabic*),ref=\arabic*]
		\item  \textbf{Functor lifting.}  For every $\C$-linear abelian category $\mathcal{C}$ and every $\C$-bilinear bifunctor
		\[
		F\colon \mathcal{A}\times\mathcal{B} \longrightarrow \mathcal{C}
		\]
		which is right exact in each variable, there exists a unique (up to natural isomorphism) $\C$-linear right exact functor
		\[
		\widehat{F} \colon \mathcal{A} \boxtimes \mathcal{B} \longrightarrow \mathcal{C}
		\]
		such that
		\[
		\widehat{F} \circ \boxtimes \; \cong \; F.
		\]
		
		\item \textbf{Natural transformation isomorphism.}   For any two $\C$-linear right exact functors
		\[
		F, G\colon \mathcal{A}\boxtimes\mathcal{B} \longrightarrow \mathcal{C},
		\]
		the induced map
		\begin{equation} \label{eq:deligne-naturality-bijection}
			\Nat(F, G)
			\longrightarrow
			\Nat(F\circ\boxtimes, \; G\circ\boxtimes),
			\qquad
			\eta \longmapsto \eta\circ\boxtimes,
		\end{equation}
		is a bijection.
	\end{enumerate}
	
	Let $R$ be a $\C$-algebra, and let $\mathfrak{M}^{f}_{R}$ denote the category of finite-dimensional left $R$-modules. The following lemma is well known; see, e.g., \cite{LF13}.
	\begin{lemma}
		Let $A$ and $B$ be finite-dimensional $\C$-algebras. Then
		\[
		\otimes \colon \mathfrak{M}^f_A \times \mathfrak{M}^f_B \longrightarrow \mathfrak{M}^f_{A \otimes B}, \quad
		(M, N) \longrightarrow M \otimes N,
		\]
		is the Deligne tensor product of $ \mathfrak{M}^f_A $ and $\mathfrak{M}^f_B $.
	\end{lemma}
	
	Recall that $\mathfrak{M}^r_V$ is the category of restricted $V$-modules.
	\begin{proposition} \label{prop:deligne-module-categories}
		Let $V$ and $U$ be $C_2$-cofinite vertex operator algebras with nonnegative integer weights.
		Then
		\[
		\otimes \colon \mathfrak{M}^r_V \times \mathfrak{M}^r_U \longrightarrow \mathfrak{M}^r_{V \otimes U}, \quad
		(M, N) \longrightarrow M \otimes N,
		\]
		is the Deligne tensor product of $ \mathfrak{M}^r_V $ and $\mathfrak{M}^r_U $.
	\end{proposition}
	
	\begin{proof}
		Let $P$ and $Q$ be projective generators of $\mathfrak{M}^r_V$ and $\mathfrak{M}^r_U$, respectively.
		Let
		\[
		A = \End_V(P)^{\opp}, \quad \text{and} \; B = \End_U(Q)^{\opp}.
		\]
		Then $A$ and $B$ are both finite-dimensional algebras.
		Moreover, the functors
		\[
		h^P = \Hom_V(P, -) \colon \mathfrak{M}^r_V \longrightarrow \mathfrak{M}^f_A,
		\]
		and
		\[
		h^Q = \Hom_U(Q, -) \colon \mathfrak{M}^r_U \longrightarrow \mathfrak{M}^f_B,
		\]
		are both equivalences of categories.
		Since $\Hom_V(P, -)$ and $\Hom_U(Q, -)$ are exact, the natural isomorphism
		\[
		\Hom_{V \otimes U}(P \otimes Q, -) \cong \Hom_V \ab\big(P, \; \Hom_U(Q, -)),
		\]
		implies that $\Hom_{V\otimes U}(P\otimes Q, -)$ is also exact.
		Hence $P\otimes Q$ is a projective generator of the category  $ \mathfrak{M}^r_{V \otimes U}$.
		
		By Lemma~\ref{lem:tensor-product-hom}, we obtain the following isomorphism of algebras:
		\[
		\End_{V \otimes U}(P \otimes Q)^{\opp}
		\cong \ab\big(\End_V(P) \otimes \End_U(Q))^{\opp}
		\cong \End_V(P)^{\opp} \otimes \End_U(Q)^{\opp}
		= A \otimes B.
		\]
		Consequently, the functor
		\[
		h^{P\otimes Q} \colon \mathfrak{M}^r_{V \otimes U} \to \mathfrak{M}^f_{A \otimes B},
		\]
		is an equivalence of categories.
		
		We obtain the following commutative diagram:
		\[
		\begin{array}{ccc}
			\mathfrak{M}^r_V \times \mathfrak{M}^r_U
			& \xrightarrow{\otimes}
			& \mathfrak{M}^r_{V \otimes U} \\[6pt]
			\Big\downarrow{\scriptstyle (h^P, h^Q)}
			& &
			\Big\downarrow{\scriptstyle h^{P \otimes Q}} \\[6pt]
			\mathfrak{M}^f_A \times \mathfrak{M}^f_B
			& \xrightarrow{\otimes}
			& \mathfrak{M}^f_{A \otimes B}
		\end{array}
		\]
		Since $(h^P, h^Q)$, $h^{P \otimes Q}$ are equivalences of categories, and since
		\[
		\otimes \colon \mathfrak{M}^f_A \times \mathfrak{M}^f_B \longrightarrow \mathfrak{M}^f_{A \otimes B}
		\]
		is the Deligne tensor product of $ \mathfrak{M}^f_A $ and $\mathfrak{M}^f_B $,
		it follows that
		\[
		\otimes \colon \mathfrak{M}^r_V \times \mathfrak{M}^r_U \longrightarrow \mathfrak{M}^r_{V \otimes U},
		\]
		is also the Deligne tensor product of $ \mathfrak{M}^r_V $ and $\mathfrak{M}^r_U $.
		The proof is complete.
		
	\end{proof}
	
	Let $D = (-)^* \colon \mathbf{Vect}_f \to (\mathbf{Vect}_f)^{\opp} $ be the dual functor.
	Let
	\[
	F \colon (\mathfrak{M}^r_V)^{\opp} \times (\mathfrak{M}^r_V)^{\opp} \longrightarrow
	\mathbf{Vect}_f
	\]
	be a bilinear bifunctor that is left exact in each variable.
	Composing $F$ with $D$, we obtain a bifunctor
	\[
	D \circ F \colon (\mathfrak{M}^r_V)^{\opp} \times (\mathfrak{M}^r_V)^{\opp}
	\longrightarrow \mathbf{Vect}_f^{\opp}.
	\]
	Equivalently, $D \circ F$  may be regarded canonically as a bifunctor from $\mathfrak{M}^r_V \times \mathfrak{M}^r_V$ to $\mathbf{Vect}_f$.
	Since $D$ is linear and exact as a contravariant functor,
	it follows that $D \circ F$,  viewed as a bifunctor on
	$ \mathfrak{M}^r_V \times \mathfrak{M}^r_V$, is bilinear and right exact in each variable.
	Hence, by the universal property of the Deligne tensor product,
	there exists a unique (up to natural isomorphism) linear right exact functor
	$\widehat{D\circ F} \colon \mathfrak{M}^r_{V \otimes V} \to \mathbf{Vect}_f$
	such that
	$ \widehat{D\circ F} \circ \otimes \; \cong \; D \circ F.$
	Equivalently, setting $ \widetilde{F} = D\circ \widehat{D\circ F}$, we obtain a unique (up to natural isomorphism) linear left exact presheaf
	$ \widetilde{F}  \colon (\mathfrak{M}^r_{V \otimes V})^{\opp} \to \mathbf{Vect}_f$
	such that
	$ \widetilde{F} \circ \otimes \cong F$.
	Let
	\[
	G \colon (\mathfrak{M}^r_V)^{\opp} \times (\mathfrak{M}^r_V)^{\opp} \longrightarrow
	\mathbf{Vect}_f
	\]
	be another bilinear bifunctor that is also left exact in each variable.
	Using the bijection in \eqref{eq:deligne-naturality-bijection} together with the natural
	isomorphism $D\circ D \cong \id$, we obtain
	\[
	\Nat(F, G) \cong \Nat \ab\big(D \circ F, \; D \circ G) \cong
	\Nat(\widehat{D\circ F}, \; \widehat{D\circ G}) \cong \Nat(\widetilde{F}, \widetilde{G}).
	\]
	
	Let $V$ be a $C_2$-cofinite vertex operator algebra with nonnegative integer weights,
	and let $W$ be a restricted $V$-module.
	Then by Proposition~\ref{prop:c2-regular-finiteness}, we have $\mathcal{R}(W) \in \mathfrak{M}^r_{V \otimes V}$.
	Since the space of intertwining operators between any three restricted $V$-modules is finite-dimensional, $\mathcal{I}^W(-, -)$ defines a bilinear bifunctor from $(\mathfrak{M}^r_V)^{\opp} \times (\mathfrak{M}^r_V)^{\opp}$ to $\mathbf{Vect}_f$, left exact in each variable.
	By the uniqueness of the lifted functor and the isomorphism \eqref{eq:hom-realization-map}, we obtain the following isomorphism of presheaves:
	\[
	\widetilde{\mathcal{I}^W} \cong h_{\mathcal{R}(W)}.
	\]
	
	\begin{theorem} \label{thm:c2-intertwiner-yoneda}
		Let $V$ be a $C_2$-cofinite vertex operator algebra with nonnegative integer weights,
		and let $W$ be a restricted $V$-module. Let
		\[
		F \colon (\mathfrak{M}^r_V)^{\opp} \times (\mathfrak{M}^r_V)^{\opp} \longrightarrow
		\mathbf{Vect}_f
		\]
		be a bilinear bifunctor that is left exact in each variable. Then there is
		a linear isomorphism
		\begin{equation} \label{eq:c2-yoneda-evaluation}
			\Nat_{(\mathfrak{M}^r_V)^{\opp}\times(\mathfrak{M}^r_V)^{\opp}}(\mathcal{I}^W, \; F) \cong \widetilde{F}(\mathcal{R}(W)).
		\end{equation}
		In particular, for every restricted $V\otimes V$-module $X$, we have the following isomorphism:
		\begin{equation} \label{eq:c2-density-hom}
			\Nat_{(\mathfrak{M}^r_V)^{\opp}\times(\mathfrak{M}^r_V)^{\opp}} \ab\big(\mathcal{H}_{\mathcal{R}(W)}, \; \mathcal{H}_X)
			\cong \Hom_{V \otimes V} \ab\big(\mathcal{R}(W), \; X).
		\end{equation}

	\end{theorem}
	
	\begin{proof}
		Together with the preceding analysis and the standard Yoneda lemma, we obtain:
		\begin{align*}
			\Nat_{(\mathfrak{M}^r_V)^{\opp}\times(\mathfrak{M}^r_V)^{\opp}}(\mathcal{I}^W, F) & \cong
			\Nat_{(\mathfrak{M}^r_V)^{\opp}\times(\mathfrak{M}^r_V)^{\opp}} \ab\big(\mathcal{H}_{\mathcal{R}(W)}, \; F) \\
			& \cong \Nat_{(\mathfrak{M}^r_{V\otimes V})^{\opp}} \ab\big(h_{\mathcal{R}(W)}, \; \widetilde{F})
			\cong \widetilde{F}(\mathcal{R}(W)).
		\end{align*}
		This completes the proof.
	\end{proof}
	
	\begin{remark}
		By \cite[Theorem~3.15]{McR23}, the Deligne tensor product (in the sense of \cite{EGNO15}) of two locally finite abelian categories of restricted modules over vertex operator algebras is realized as a certain subcategory of the category of restricted modules over the tensor product vertex operator algebra.
		We point out, however, that $\mathcal{R}(W)$ is usually not contained in the Deligne tensor product.
	\end{remark}
	
	\section{A Yoneda-Type Lemma and Coend Realization}\label{sec:yoneda-and-coend}
	
	The first main objective of this chapter is to extend the isomorphism \eqref{eq:c2-density-hom} from $C_2$-cofinite restricted modules to a general framework.
	This extension is established in Lemma \ref{lem:external-tensor-density}.
	The resulting generalized isomorphism has three direct applications, which are presented in
	Theorem~\ref{thm:intertwiner-yoneda}, Theorem~\ref{thm:regular-coend-realization}, and Theorem~\ref{thm:li-regular-comparison}.
	
	The proof of the isomorphism \eqref{eq:c2-density-hom} using the theory of Deligne tensor products involves the categories $\mathfrak{M}^r_V \times \mathfrak{M}^r_V$ and $\mathfrak{M}^r_{V\otimes V}$, together with the $V\otimes V$-module $\mathcal{R}(W)$. We generalize this framework by replacing $\mathfrak{M}^r_V \times \mathfrak{M}^r_V$ with the product category $\mathcal{A}\times\mathcal{B}$, $\mathfrak{M}^r_{V\otimes V}$ with the abelian category $\mathcal{G}_{\mathcal{A}, \mathcal{B}}$ generated by $\mathcal{A}\odot\mathcal{B}$, and $\mathcal{R}(W)$ with the distinguished submodule $\mathcal{R}_{\mathcal{A}, \mathcal{B}}(W) \subseteq \mathcal{R}(W)$, where $\mathcal{A}$ and $\mathcal{B}$ are full subcategories of the category of weak $V$-modules.
	These notions are introduced in Section~\ref{subsec:generated-modules}.
	We note that the theory of Deligne tensor products is generally not available in this extended setting.
	
	Although our actual application concerns $V\otimes V$-modules,
	we work with general $V\otimes U$-modules. This generalization is arguably more natural and adds no extra complexity to our arguments.
	
	\subsection{ The Category \texorpdfstring{$ \mathcal{G}_{\mathcal{A}, \mathcal{B}}$}{G\_(A, B)} and The Distinguished Submodule  \texorpdfstring{$\mathcal{S}_{\mathcal{A}, \mathcal{B}}(X) $}{S\_(A, B)(X)} } \label{subsec:generated-modules}

	Let $V$ be a vertex operator algebra, and let $\mathcal{A}$ be a full subcategory of the category of weak $V$-modules.
	We say that $\mathcal{A}$ is \textbf{closed under submodules} (resp.\ \textbf{closed under quotients}) if every $V$-submodule (resp.\ quotient module) of an object of $\mathcal{A}$ is again an object of $\mathcal{A}$. Similarly, we say that $\mathcal{A}$ is \textbf{closed under finite direct sums} (resp.\ \textbf{closed under arbitrary direct sums}) if every finite (resp.\ arbitrary) direct sum of objects of $\mathcal{A}$ is again an object of $\mathcal{A}$.
	Assume that the full subcategory $\mathcal{A}$ lies inside the category of restricted $V$-modules.
	We say that $\mathcal{A}$ is \textbf{closed under restricted duals} if every restricted dual of an object of $\mathcal{A}$ is again an object of $\mathcal{A}$.

	\vspace{0.2cm}
	
	Let $V$ and $U$ be vertex operator algebras.
	Let $\mathcal{A}$ (resp.\ $\mathcal{B}$) be a nonempty full subcategory of the category of weak $V$-modules (resp.\ weak $U$-modules).
	Consider the product category $\mathcal{A} \times \mathcal{B}$.
	Define the class
	\[
	\mathcal{A} \odot \mathcal{B} = \ab\{ M\otimes N \mid (M, N) \in \mathcal{A}\times\mathcal{B} \},
	\]
	and let $\mathcal{G}_{\mathcal{A}, \mathcal{B}}$  denote the full subcategory of weak $V\otimes U$-modules generated by $\mathcal{A} \odot \mathcal{B}$.
	Equivalently, a weak $V\otimes U$-module $X$ lies in $\mathcal{G}_{\mathcal{A}, \mathcal{B}}$  if and only if there exist an indexed family $\{M_i \otimes N_i\}_{i \in I}$ of objects in $\mathcal{A} \odot \mathcal{B}$ and a surjective homomorphism of weak $V\otimes U$-modules
	\[
	\bigoplus_{i \in I} (M_i\otimes N_i) \twoheadrightarrow X.
	\]
	Clearly,  the full subcategory $\mathcal{G}_{\mathcal{A}, \mathcal{B}}$ contains the zero module and is closed under arbitrary direct sums and quotients.
	
	
	For a weak $V\otimes U$-module $X$, define
	\[
	\mathcal{S}_{\mathcal{A}, \mathcal{B}}(X)
	= \sum \ab\big\{ \Image f \mid f \in \Hom_{V\otimes U}(M\otimes N, \; X), \ M\otimes N \in \mathcal{A} \odot \mathcal{B} \}
	\subseteq X.
	\]
	Then $\mathcal{S}_{\mathcal{A}, \mathcal{B}}(X)$ is the unique largest submodule of $X$ generated by $\mathcal{A} \odot \mathcal{B}$.
	Consequently, $X \in \mathcal{G}_{\mathcal{A}, \mathcal{B}}$ if and only if $\mathcal{S}_{\mathcal{A}, \mathcal{B}}(X) = X$.
	Moreover, if $X$ is a weak $V\otimes U$-module and $Z \in \mathcal{G}_{\mathcal{A}, \mathcal{B}}$, then the induced map
	\begin{equation}
		\Hom_{V \otimes U} \ab\big(Z, \; \mathcal{S}_{\mathcal{A}, \mathcal{B}}(X)) \hookrightarrow \Hom_{V \otimes U} (Z, \; X)
	\end{equation}
	is a linear isomorphism.
	When $U = V$, for notational simplicity, we write
	\[
	\mathcal{R}_{\mathcal{A}, \mathcal{B}}(W)
	:=
	\mathcal{S}_{\mathcal{A}, \mathcal{B}} \ab\big(\mathcal{R}(W)).
	\]

	\begin{lemma} \label{lem:pure-tensor-presentation}
		Let $V, U$, and $\mathcal{A}, \mathcal{B}$ be as above.
		Suppose that $\mathcal{A}$ and $\mathcal{B}$ are closed under finite direct sums.
		Let $X$ be a weak $V\otimes U$-module.
		Then $X \in \mathcal{G}_{\mathcal{A}, \mathcal{B}}$ if and only if for every $m \in X$,
		there exist an object $M \otimes N \in \mathcal{A} \odot \mathcal{B}$ (depending on $m$) , an element $ a \otimes b \in M \otimes N$, and a $V \otimes U$-module homomorphism
		\[
		f\colon M \otimes N \longrightarrow X
		\]
		such that $f(a \otimes b) = m$.
		
	\end{lemma}
	
	\begin{proof}
		
		($\Leftarrow$ ) Since  $\mathcal{G}_{\mathcal{A}, \mathcal{B}}$ is closed under arbitrary direct sums and quotients,  the assertion in this direction follows immediately.
		
		~~~~~( $ \Rightarrow $ ) Since $X \in \mathcal{G}_{\mathcal{A}, \mathcal{B}}$, for every $m \in X$, there exist finitely many $V \otimes U$-modules
		\[
		M_1 \otimes N_1, \dots, M_n \otimes N_n \in \mathcal{A} \odot \mathcal{B},
		\]
		and a $V \otimes U$-module homomorphism
		\[
		g \colon \bigoplus_{i = 1}^n M_i \otimes N_i \to X
		\]
		such that
		\[
		g \ab\Big(\sum_{i = 1}^n \sum_{j = 1}^{k_i} a_{i, j} \otimes b_{i, j}) = m.
		\]
		for some $a_{i, j} \in M_i, b_{i, j} \in N_i$, and some positive integers $k_i$.
		Set
		\[
		M := \bigoplus_{i = 1}^n M_i^{\oplus k_i}, \;\;\;\; N := \bigoplus_{i = 1}^n N_i^{ \oplus k_i},
		\]
		and let
		\[
		a := (a_{1, 1}, \dots, a_{n, k_n}) \in M, \; \;\;\; b := (b_{1, 1}, \dots, b_{n, k_n}) \in N.
		\]
		By our assumptions on $\mathcal{A}$ and $\mathcal{B}$, we have $M \in \mathcal{A}$ and $N \in \mathcal{B}$.
		Define a linear map
		\[
		f \colon M \otimes N \longrightarrow X, \;\; \;\;
		\]
		by setting
		\[
		f \ab\Big( \ab\big(\widetilde{a}_{1, 1}, \dots, \widetilde{a}_{i, j}, \dots, \widetilde{a}_{n, k_n}) \otimes
		\ab\big(\widetilde{b}_{1, 1}, \dots, \widetilde{b}_{i, j}, \dots, \widetilde{b}_{n, k_n}) ) = \sum_{i = 1}^n \sum_{j = 1}^{k_i} g(\widetilde{a}_{i, j} \otimes \widetilde{b}_{i, j}).
		\]
		for all $(\widetilde{a}_{i, j}) \in M$ and $(\widetilde{b}_{i, j}) \in N$.  Then $f$ is a $V \otimes U$-module homomorphism, and
		\[
		f(a \otimes b) = g \ab\Big(\sum_{i = 1}^n \sum_{j = 1}^{k_i} a_{i, j} \otimes b_{i, j}) = \sum_{i = 1}^n \sum_{j = 1}^{k_i} g \ab\big(a_{i, j} \otimes b_{i, j}) = m.
		\]
		This proves the direction ($ \Rightarrow $).
	\end{proof}
	
	\begin{lemma} \label{lem:common-presentation-refinement}
		Let $V, U$, and $\mathcal{A}, \mathcal{B}$ be as above.
		Suppose that $\mathcal{A}$ and $\mathcal{B}$ are closed under finite direct sums and submodules.
		Let $X$ be a weak $V \otimes U$-module, and let $m \in X$.
		Let $(M_0, N_0)$ and $(M_1, N_1)$ be objects of $ \mathcal{A} \times \mathcal{B}$. Suppose that
		\[
		f_0 \colon M_0 \otimes N_0 \longrightarrow X, \qquad f_1 \colon M_1 \otimes N_1 \longrightarrow X
		\]
		are $V \otimes U$-module homomorphisms satisfying
		\[
		f_0(a_0 \otimes b_0) = f_1(a_1 \otimes b_1) = m
		\]
		for some $a_0 \otimes b_0 \in M_0 \otimes N_0$ and $a_1 \otimes b_1 \in M_1 \otimes N_1$.
		
		Then there exist an object $(M, N) \in \mathcal{A} \times \mathcal{B} $, an element $a \otimes b \in M \otimes N$,
		and two morphisms
		\[
		(\alpha_0, \beta_0) \colon (M, N) \longrightarrow (M_0, N_0), \qquad
		(\alpha_1, \beta_1) \colon (M, N) \longrightarrow (M_1, N_1)
		\]
		in $\mathcal{A} \times \mathcal{B} $ such that
		\[
		(\alpha_0 \otimes \beta_0) (a \otimes b) = a_0 \otimes b_0, \qquad
		(\alpha_1 \otimes \beta_1) (a \otimes b) = a_1 \otimes b_1.
		\]
		Moreover, the following diagram is commutative:
		\[
		\begin{tikzcd}
			M \otimes N \arrow[r, "\alpha_0 \otimes \beta_0"] \arrow[d, "\alpha_1 \otimes \beta_1"'] & M_0 \otimes N_0 \arrow[d, "f_0"] \\
			M_1 \otimes N_1 \arrow[r, "f_1"'] & X
		\end{tikzcd}
		\]
		
	\end{lemma}
	
	\begin{proof}
		
		Let $M_2 = M_0 \oplus M_1$ and $N_2 = N_0 \oplus N_1$, and set $a = (a_0, a_1) \in M_2$ and $b = (b_0, b_1) \in N_2$.
		Let $M$ be the cyclic $V$-submodule of $M_2$ generated by $(a_0, a_1)$, and let  $N$ be the cyclic $U$-submodule of $N_2$ generated by $(b_0, b_1)$.
		Then $M \otimes N$ is the cyclic $V \otimes U$-submodule of $M_2 \otimes N_2$ generated by $a \otimes b$.
		Moreover, since $\mathcal{A}$ and $\mathcal{B}$ are closed under submodules, we have $(M, N) \in \mathcal{A} \times \mathcal{B} $.
		Let
		\[
		p_0 \colon M_2 \to M_0, \quad p_1 \colon M_2 \to M_1, \quad
		q_0 \colon N_2 \to N_0, \quad q_1 \colon N_2 \to N_1
		\]
		be the canonical projections onto the respective direct-summand factors.
		Denote by $\iota_M\colon M \hookrightarrow M_2$ and $\iota_N\colon N \hookrightarrow N_2$ the respective inclusion homomorphisms, and define
		\[
		\alpha_0 = p_0 \circ \iota_M, \quad \beta_0 = q_0 \circ \iota_N, \quad \alpha_1 = p_1 \circ \iota_M, \quad \beta_1 = q_1 \circ \iota_N.
		\]
		A straightforward verification shows that
		\[
		\ab\big(\alpha_0 \otimes \beta_0) (a \otimes b) = a_0 \otimes b_0, \qquad \ab\big(\alpha_1 \otimes \beta_1) (a \otimes b) = a_1 \otimes b_1.
		\]
		Since $M\otimes N$ is generated as a $V\otimes U$-module by the single element $a\otimes b$, the equality
		\[
		f_0 \circ \ab\big(\alpha_0 \otimes \beta_0) (m_0) = f_1 \circ \ab\big(\alpha_1 \otimes \beta_1) (m_0)
		\]
		holds for all $m_0 \in M \otimes N$.
		This is exactly the condition that the diagram above commutes.
		
	\end{proof}
	
	\begin{corollary} \label{cor:generated-category-closure}
		
		Let $V, U$, and $\mathcal{A}, \mathcal{B}$ be as above.
		Suppose that $\mathcal{A}$ and $\mathcal{B}$ are closed under finite direct sums and submodules.
		Then the category $\mathcal{G}_{\mathcal{A}, \mathcal{B}}$ is closed under arbitrary direct sums, submodules, and quotients.
		
	\end{corollary}
	
	\begin{proof}
		
		The closure of $\mathcal{G}_{\mathcal{A}, \mathcal{B}}$ under arbitrary direct sums and quotients follows immediately from its definition.
		Let $X \in \mathcal{G}_{\mathcal{A}, \mathcal{B}}$, and let $Z$ be a submodule of $X$.
		By Lemma \ref{lem:pure-tensor-presentation}, every element $m \in Z$ is in the image of some homomorphism $X_0 \to Z$ with
		$X_0 \in \mathcal{A}\odot\mathcal{B}$.
		It follows that $Z$ is generated by $\mathcal{A}\odot\mathcal{B}$, so that $Z \in \mathcal{G}_{\mathcal{A}, \mathcal{B}}$.
		This shows that $\mathcal{G}_{\mathcal{A}, \mathcal{B}}$ is closed under taking submodules.
		The proof is complete.
		
	\end{proof}
	
	\begin{remark}
		It follows from Corollary~\ref{cor:generated-category-closure} that $\mathcal{G}_{\mathcal{A}, \mathcal{B}}$
		is in fact an exact abelian subcategory of the category of weak $V\otimes U$-modules.
		In other words, $\mathcal{G}_{\mathcal{A}, \mathcal{B}}$ is a full subcategory which is abelian, and the inclusion functor is exact.
	\end{remark}

	\begin{definition}
		Let $V$ and $U$ be vertex operator algebras.
		A generalized $V$-module is called \textbf{locally restricted} if every finitely generated submodule is a restricted module.
		A locally restricted $V$-module is said to be of \textbf{locally finite length} if every finitely generated submodule has finite composition length.
		A weak $V\otimes U$-module $X$ is called a \textbf{bi-generalized $V\otimes U$-module} if $X$ is generalized both as a $V$-module and as a $U$-module.
		In particular, $X$ is a generalized $V\otimes U$-module. 
	\end{definition}

	\begin{proposition} \label{prop:generated-category-characterizations}
		Let $V$ and $U$ be vertex operator algebras, and let $X$ be a weak $V\otimes U$-module. The following statements hold.
		\begin{enumerate}[label=(\arabic*),ref=\arabic*]
			\item Let $\mathcal{A}$ (resp.\ $\mathcal{B}$) denote the category of weak $V$-modules (resp.\ weak $U$-modules).
			Then $\mathcal{G}_{\mathcal{A}, \mathcal{B}}$ is precisely the category of weak $V\otimes U$-modules, and $\mathcal{S}_{\mathcal{A}, \mathcal{B}}(X) = X$.
			
			\item Let $\mathcal{A}$ (resp.\ $\mathcal{B}$) denote the category of generalized $V$-modules (resp.\ generalized $U$-modules).
			Then $\mathcal{G}_{\mathcal{A}, \mathcal{B}}$ is precisely the category of bi-generalized $V\otimes U$-modules,
			and $\mathcal{S}_{\mathcal{A}, \mathcal{B}}(X)$ is the unique largest bi-generalized $V\otimes U$-submodule of $X$.
			
			\item Let $\mathcal{A}$ (resp.\ $\mathcal{B}$) denote the category of restricted $V$-modules (resp.\ restricted $U$-modules).
			Then $\mathcal{G}_{\mathcal{A}, \mathcal{B}}$ is precisely the category of locally restricted $V\otimes U$-modules,
			and $\mathcal{S}_{\mathcal{A}, \mathcal{B}}(X)$ is the sum of all restricted $V\otimes U$-submodules of $X$.
			
		\end{enumerate}
	\end{proposition}
	
	\begin{proof}
		
		Since $\mathcal{S}_{\mathcal{A}, \mathcal{B}}(X)$ is, by definition, the unique largest submodule of $X$ generated by $\mathcal{A}\odot\mathcal{B}$,
		it is automatically the unique largest submodule of $X$ belonging to $\mathcal{G}_{\mathcal{A}, \mathcal{B}}$.
		Consequently, once $\mathcal{G}_{\mathcal{A}, \mathcal{B}}$ has been characterized,
		the description of $\mathcal{S}_{\mathcal{A}, \mathcal{B}}(X)$ follows immediately.
		
		\vspace{0.2cm}
		
		(1) Let $X$ be a weak $V \otimes U$-module. By Lemma~\ref{lem:cyclic-tensor-presentation}, every element $m \in X$ is in the image of some homomorphism $Z \to X$ with $Z \in \mathcal{A}\odot\mathcal{B}$.
		It follows that $X \in \mathcal{G}_{\mathcal{A}, \mathcal{B}}$. The result follows.
		
		\vspace{0.2cm}
		
		(2) Since each object $M\otimes N \in \mathcal{A}\odot\mathcal{B}$ is a bi-generalized $V\otimes U$-module, and arbitrary direct sums and quotients of bi-generalized $V\otimes U$-modules are again bi-generalized, it follows that every object of $\mathcal{G}_{\mathcal{A}, \mathcal{B}}$ is bi-generalized.
		
		Let $X$ be a bi-generalized $V\otimes U$-module. By Lemma~\ref{lem:cyclic-tensor-presentation}, for any $m \in X$, there exist a $V$-submodule $M$ of $X$, a $U$-submodule $N$ of $X$, elements $a \in M$, $b \in N$, and a $V\otimes U$-module homomorphism $f:M\otimes N \to X$ such that $f(a\otimes b) = m$.
		Since $M$ is a $V$-submodule of $X$, we have $M \in \mathcal{A}$, and similarly $N \in \mathcal{B}$.
		Consequently, every element $m \in X$ is in the image of some homomorphism $Z \to X$ with $Z \in \mathcal{A}\odot\mathcal{B}$.
		It follows that $X$ is generated by
		$\mathcal{A}\odot \mathcal{B}$.
		
		\vspace{0.2cm}
		
		(3) We first observe that for any $(M, N) \in \mathcal{A}\times\mathcal{B}$, $M\otimes N$ is a locally restricted $V\otimes U$-module.
		Note, however, that $M$ or $N$ need not be finitely generated, so $M\otimes N$ may fail to be a restricted $V\otimes U$-module.
		Since the property of being locally restricted is closed under arbitrary direct sums and quotients, every object in $\mathcal{G}_{\mathcal{A}, \mathcal{B}}$ is locally restricted.
		
		Let $X$ be a locally restricted $V\otimes U$-module.
		Since $L(0)$ acts locally finitely on $X$, the operators $L^V(0)$ and $L^U(0)$ also act locally finitely on $X$.
		Consequently, $X$ is a bi-generalized $V \otimes U$-module.
		
		Let $m$ be an arbitrary element of $X$. Without loss of generality,
		we may assume that $m$ is a common generalized eigenvector of $L^V(0)$ and $L^U(0)$.
		Let $Z$ denote the cyclic $V\otimes U$-module generated by $m$.
		Let $M$ (resp. $N$) stand for the cyclic $V$-submodule (resp. $U$-submodule) of $X$ generated by $m$.
		Since $Z$ is a restricted $V\otimes U$-module, and $M$ (resp.\ $N$) is a generalized $V$-module (resp.\ $U$-module);
		hence both $M$ and $N$ are necessarily restricted modules.
		By Lemma~\ref{lem:cyclic-tensor-presentation}, $Z$ is a quotient of $M\otimes N$, and so $Z \in \mathcal{G}_{\mathcal{A}, \mathcal{B}}$.
		Since $\mathcal{G}_{\mathcal{A}, \mathcal{B}}$ is closed under arbitrary direct sums and quotients.
		Therefore, $X \in \mathcal{G}_{\mathcal{A}, \mathcal{B}}$.
		This completes the proof.
		
	\end{proof}
	
	\begin{lemma} \label{lem:c1-cofinite-cyclic-factors}
		
		Let $V$ (resp.\ $U$) be an $\N$-graded vertex operator algebra such that every submodule of a $C_1$-cofinite restricted $V$-module (resp.\ $C_1$-cofinite restricted $U$-module) is again $C_1$-cofinite.
		Let $X$ be a $C_1$-cofinite restricted $V\otimes U$-module.
		Let $m \in X$, and let $M$ (resp.\ $N$) denote the cyclic $V$-submodule (resp.\ $U$-submodule) of $X$ generated by $m$.
		Then $M$ (resp.\ $N$) is a $C_1$-cofinite restricted $V$-module (resp.\ $U$-module).
		
	\end{lemma}
	
	\begin{proof}
		This is established in the second through penultimate paragraphs of the proof of \cite[Theorem 3.7]{McR23}.
		
	\end{proof}
	
	\begin{proposition} \label{prop:generated-locally-finite-length-modules}
		Let $V$ (resp.\ $U$) be an $\N$-graded vertex operator algebra such that
		a restricted $V$-module (resp.\ a restricted $U$-module) is $C_1$-cofinite
		if and only if it has finite composition length.
		Let $\mathcal{A}$ (resp.\ $\mathcal{B}$) denote the category of restricted $V$-modules (resp.\ restricted $U$-modules) of finite composition length.
		Then $\mathcal{G}_{\mathcal{A}, \mathcal{B}}$ is the category of locally restricted $V\otimes U$-modules of locally finite length, and  for every weak $V \otimes V$-module $X$, $\mathcal{S}_{\mathcal{A}, \mathcal{B}}(X)$ is precisely the sum of all restricted $V\otimes U$-submodules of $X$ of finite composition length.
	\end{proposition}
	
	\begin{proof}
		
		We first observe that for any $(M, N) \in \mathcal{A}\times\mathcal{B}$, the tensor product $M\otimes N$ is a $V\otimes U$-module of finite length.
		Consequently, every object in $ \mathcal{G}_{\mathcal{A}, \mathcal{B}}$ has locally finite length.
		To prove that every locally restricted $V\otimes U$-modules of locally finite length belongs to $ \mathcal{G}_{\mathcal{A}, \mathcal{B}}$,
		it suffices to show that all restricted $V\otimes U$-modules of finite composition length lie in $ \mathcal{G}_{\mathcal{A}, \mathcal{B}}$.
		
		Let $X$ be a restricted $V\otimes U$-module of finite composition length.
		Then $X$ is $C_1$-cofinite. Indeed, every simple restricted $V\otimes U$-module is isomorphic to $M\otimes N$, where $M$ is a simple restricted $V$-module and $N$ is a simple restricted $U$-module\cite{FHL93}. By our assumptions on $V$ and $U$, $M$ and $N$ are both $C_1$-cofinite. Since the tensor product of two $C_1$-cofinite modules remains $C_1$-cofinite, every simple restricted $V\otimes U$-module is $C_1$-cofinite. Recall that $C_1$-cofiniteness is preserved under finite extensions. Consequently, any restricted $V\otimes U$-module of finite composition length is $C_1$-cofinite. In particular, $X$ is $C_1$-cofinite.
		
		Let $m \in X$, and let $M$ (resp. $N$) denote the cyclic $V$-submodule (resp. $U$-submodule) of $X$ generated by $m$.
		By Lemma~\ref{lem:c1-cofinite-cyclic-factors}, $M$ and $N$ are both $C_1$-cofinite and hence of finite composition length.
		By Lemma~\ref{lem:cyclic-tensor-presentation},
		$m$ lies in the image of a homomorphism $M\otimes N \to X$. It follows that $X$ belongs to $\mathcal{G}_{\mathcal{A}, \mathcal{B}}$. This completes the proof.
		
	\end{proof}
	
	\begin{remark}
		Every $C_2$-cofinite vertex operator algebra without negative conformal weights satisfies the hypotheses of
		Property~\ref{prop:generated-locally-finite-length-modules}.
		In addition, it is conjectured in \cite[Conjecture~3.3]{CY21} that
		if $V$ is a vertex operator algebra such that every irreducible ordinary $V$-module is $C_1$-cofinite,
		then a restricted $V$-module is $C_1$-cofinite if and only if it has finite length.
		Moreover, if the vertex operator algebra $V$ satisfies the assumptions of Property~\ref{prop:generated-locally-finite-length-modules},
		then the category of restricted $V$-modules of finite composition length forms a braided tensor category \cite{CMHFY26}.
		
	\end{remark}
	
	\subsection{ A Yoneda-Type Lemma }
	
	Let $V$ and $U$ be vertex operator algebras.
	Let $\mathcal{A}$ (resp.\ $\mathcal{B}$) be a nonempty full subcategory of the category of weak $V$-modules (resp.\ weak $U$-modules).
	Let $X$ be a weak $V\otimes U$-module.
	Then $\mathcal{H}_X=\Hom_{V\otimes U}(-\otimes -, X)$ is a bifunctor, contravariant in each module variable, from $\mathcal{A}^{\opp}\times\mathcal{B}^{\opp}$ to $\mathbf{Vect}$.
	Let $Z$ be a weak $V\otimes U$-module, and let $\beta\colon X \to Z$ be a $V\otimes U$-module homomorphism.
	Then $\beta$ induces a natural transformation
	\[
	\mathcal{H}(\beta) \colon \mathcal{H}_X \Rightarrow \mathcal{H}_Z.
	\]
	Its component at $(M,N)$ sends $f$ to $\beta\circ f$.
	The following lemma can be regarded as a weak form of the Yoneda lemma.
	
	\begin{lemma} \label{lem:external-tensor-density}
		
		Let $U, V$, $\mathcal{A}$, and $\mathcal{B}$ be as above.
		Suppose that $\mathcal{A}$ and $\mathcal{B}$ are closed under finite direct sums and submodules.
		Let $X \in \mathcal{G}_{\mathcal{A}, \mathcal{B}}$, and let $Z$ be a weak $V\otimes U$-module.
		Then the map
		\[
		\Theta\colon \Hom_{V\otimes U}(X, \; Z) \longrightarrow
		\Nat_{\mathcal{A}^{\opp}\times\mathcal{B}^{\opp}} \ab\big(\mathcal{H}_X, \; \mathcal{H}_Z)
		\]
		defined by
		\[
		\Theta(\beta) = \mathcal{H}(\beta)
		\]
		is a linear isomorphism.
		Moreover, if we further assume $Z \in \mathcal{G}_{\mathcal{A}, \mathcal{B}}$,
		then $\Theta(\beta)$ is a natural isomorphism of bifunctors if and only if $\beta$ is a $V \otimes U$-module isomorphism.  
	\end{lemma}
	
	\begin{proof}
		
		We divide the proof into three parts: injectivity of $\Theta$, surjectivity of $\Theta$, and the equivalence for isomorphisms.
		
		\smallskip
		\noindent\textbf{Part 1: $\Theta$ is injective.}
		Suppose $\Theta(\beta) = 0$ for some $\beta \in \Hom_{V \otimes U} (X, Z)$.
		Then for every $M \otimes N \in \mathcal{A} \odot \mathcal{B}$, and any $V \otimes U$-module homomorphism $f \colon M \otimes N \to X$, we have $\beta f = 0$.
		Since $X$ is gnerated by $\mathcal{A} \odot \mathcal{B}$, this forces $\beta = 0$.
		Therefore $\Theta$ is injective.  Linearity of $\Theta$ is immediate from the definition.
		
		\smallskip
		\noindent\textbf{Part 2: $\Theta$ is surjective.}
		Let
		\[
		\Gamma \colon \mathcal{H}_X \Rightarrow \mathcal{H}_Z
		\]
		be an arbitrary natural transformation.
		Throughout the proof below, $\Gamma_{M, N}$ denotes the component of $\Gamma$ at $(M, N)$.
		Our goal is to construct a $V\otimes U$-module homomorphism $\beta_{\Gamma}\colon X \to Z$ such that
		\[
		\Gamma = \mathcal{H}(\beta_\Gamma).
		\]
		
		Let $m \in X$.  By Lemma~\ref{lem:pure-tensor-presentation}, there exists an object $M \otimes N$ of $\mathcal{A} \odot \mathcal{B}$,
		an element $a \otimes b \in M \otimes N$, and a $V \otimes U$-module homomorphism $f \colon M \otimes N \to X$ such that $f(a \otimes b) = m$.
		We define
		\[
		\beta_\Gamma(m) = \Gamma_{M, N}(f) (a \otimes b).
		\]
		
		We first show that this definition is independent of the choice of the triple $(M\otimes N, \; a \otimes b, \; f)$.
		Take two such triples $(M_i\otimes N_i, \; a_i\otimes b_i, \; f_i)$, $i = 0, 1$, satisfying $f_i(a_i\otimes b_i) = m$.
		By Lemma \ref{lem:common-presentation-refinement}, there exists an object $M \otimes N \in \mathcal{A} \odot \mathcal{B}$,
		an element $a \otimes b \in M \otimes N$,
		and two morphisms
		\[
		(g_0, h_0) \colon (M, N) \longrightarrow (M_0, N_0), \qquad
		(g_1, h_1) \colon (M, N) \longrightarrow (M_1, N_1)
		\]
		in $\mathcal{A} \times \mathcal{B} $ such that
		\[
		(g_0 \otimes h_0) (a \otimes b) = a_0 \otimes b_0, \;\; (g_1 \otimes h_1) (a \otimes b) = a_1 \otimes b_1,
		\;\; \text{and} \;\; f_0 (g_0 \otimes h_0) = f_1 (g_1 \otimes h_1).
		\]
		Applying the naturality of $\Gamma$ to the morphism $g_0 \otimes h_0 \colon M \otimes N \longrightarrow M_0 \otimes N_0$, we obtain the following commutative diagram:
		\[
		\begin{tikzcd}
			\Hom_{V \otimes U}(M_0 \otimes N_0, X)
			\arrow[r, "{\Gamma_{M_0, N_0}}"]
			\arrow[d, "{\Hom_{V \otimes U}(g_0 \otimes h_0, X)}"']
			&
			\Hom_{V \otimes U}(M_0 \otimes N_0, Z)
			\arrow[d, "{\Hom_{V \otimes U}(g_0 \otimes h_0, Z)}"] \\
			\Hom_{V \otimes U}(M \otimes N, X)
			\arrow[r, "{\Gamma_{M, N}}"']
			&
			\Hom_{V \otimes U}(M \otimes N, Z)
		\end{tikzcd}
		\]
		It follows from the commutative diagram that
		\begin{align*}
			\Gamma_{M_0, N_0}(f_0) (a_0 \otimes b_0) & = \Gamma_{M_0, N_0}(f_0) (g_0 \otimes h_0) (a \otimes b) \\
			& = \ab\big(\ab\big(\Hom_{V \otimes U}(g_0 \otimes h_0, Z) \circ \Gamma_{M_0, N_0}) (f_0)) (a \otimes b) \\
			& = \ab\big(\ab\big(\Gamma_{M, N} \circ \Hom_{V \otimes U}(g_0 \otimes h_0, X)) (f_0)) (a \otimes b) \\
			& = \Gamma_{M, N} \ab\big(f_0 (g_0 \otimes h_0)) (a \otimes b).
		\end{align*}
		Similarly, applying the naturality of $\Gamma$ to the morphism $g_1 \otimes h_1\colon M \otimes N \to M_1\otimes N_1$ gives
		\[
		\Gamma_{M_1, N_1}(f_1) (a_1 \otimes b_1) = \Gamma_{M, N} \ab\big(f_1 (g_1 \otimes h_1)) (a \otimes b).
		\]
		Since $f_0 (g_0 \otimes h_0) = f_1 (g_1 \otimes h_1)$, we obtain
		\[
		\Gamma_{M_0, N_0}(f_0) (a_0 \otimes b_0) = \Gamma_{M_1, N_1}(f_1) (a_1 \otimes b_1).
		\]
		Hence $\beta_{\Gamma}$ is well-defined.
		
		\vspace{0.2cm}
		
		By definition of $\beta_{\Gamma}$, for every $V \otimes U$-module homomorphism $f \colon M \otimes N \to X$, and every $a \otimes b \in M \otimes N$
		we have
		\begin{align} \label{eq:density-component-evaluation}
			\Gamma_{M, N}(f) (a \otimes b) = \beta_{\Gamma} \ab\big(f (a \otimes b)).
		\end{align}
		Since this equality holds for all elements $a \otimes b \in M \otimes N$, we conclude
		\[
		\Gamma_{M, N}(f) = \beta_{\Gamma} \circ f,
		\]
		for all $V \otimes U$-module homomorphism $f \colon M \otimes N \to X$.
		
		\vspace{0.2cm}
		
		Next we show that $\beta_{\Gamma} $  is linear.
		Let $m_0, m_1 \in X$.
		For $i = 0, 1$, choose an object $M_i \otimes N_i$ of $\mathcal{A} \odot \mathcal{B}$,
		an element $a_i \otimes b_i \in M_i \otimes N_i$, and a $V \otimes U$-module homomorphism $f_i \colon M_i \otimes N_i \to X$ such that
		$f_i(a_i \otimes b_i) = m_i$.
		Let
		\[
		M := M_0 \oplus M_1, \quad N := N_0 \oplus N_1,
		\]
		and define a $V \otimes U$-homorphism
		\[
		f \colon M \otimes N \longrightarrow X, \quad (\widetilde{a}_0, \widetilde{a}_1) \otimes (\widetilde{b}_0, \widetilde{b}_1) \longrightarrow f_0(\widetilde{a}_0 \otimes \widetilde{b}_0) + f_1(\widetilde{a}_1 \otimes \widetilde{b}_1).
		\]
		Then we have
		\begin{align*}
			& f \ab\big((a_0, 0) \otimes (b_0, 0)) = f_0(a_0 \otimes b_0) = m_0,
			\quad f \ab\big((0, a_1) \otimes (0, b_1)) = f_1(a_1 \otimes b_1) = m_1; \\
			& f \ab\big((a_0, 0) \otimes (0, b_1)) = f \ab\big((0, a_1) \otimes (b_0, 0)) = 0;\\
			& f \ab\big((a_0, a_1) \otimes (b_0, b_1)) = m_0 + m_1.
		\end{align*}
		Using the linearity of $\Gamma_{M, N}(f)$ and Identity \eqref{eq:density-component-evaluation},
		\begin{align*}
			\beta_{\Gamma} \ab\big(m_0 + m_1) & = \beta_{\Gamma} \ab\big(f \ab\big((a_0, a_1) \otimes (b_0, b_1))) \\
			& = \Gamma_{M, N}(f) \ab\big((a_0, a_1) \otimes (b_0, b_1)) \\
			& = \Gamma_{M, N}(f) \ab\big((a_0, 0) \otimes (b_0, 0)) + \Gamma_{M, N}(f) \ab\big((0, a_1) \otimes (0, b_1)) \\
			& \quad +\Gamma_{M, N}(f) \ab\big((a_0, 0) \otimes (0, b_1)) +\Gamma_{M, N}(f) \ab\big((0, a_1) \otimes (b_0, 0)) \\
			& = \beta_{\Gamma} \ab\big(f \ab\big((a_0, 0) \otimes (b_0, 0)))
			+ \beta_{\Gamma} \ab\big(f \ab\big((0, a_1) \otimes (0, b_1))) \\
			& \quad +\beta_{\Gamma} \ab\big(f \ab\big((a_0, 0) \otimes (0, b_1))) + \beta_{\Gamma} \ab\big(f \ab\big((0, a_1) \otimes (b_0, 0))) \\
			& = \beta_{\Gamma}(m_0) +\beta_{\Gamma} (m_1).
		\end{align*}
		This establishes that $\beta_\Gamma$ preserves addition. Preservation of scalar multiplication is straightforward, so $\beta_\Gamma$ is linear.
		
		\vspace{0.2cm}
		
		It remains to check that $\beta_\Gamma$ is a $V\otimes U$-module homomorphism.
		Let $m \in X$,
		and choose an object $M \otimes N \in \mathcal{A} \odot \mathcal{B}$, an element $a \otimes b \in M \otimes N$, and a $V \otimes U$-module homomorphism
		$f \colon M \otimes N \to X$ such that $f(a \otimes b) = m$.
		For all $v \in V$ and $ u \in U$, we have
		\begin{align*}
			\beta_{\Gamma} \ab\big((v \otimes u)_n m) & = \beta_{\Gamma} \ab\big((v \otimes u)_n f(a \otimes b))
			= \beta_{\Gamma} \ab\bigg(f \ab\big((v \otimes u)_n (a \otimes b))) \\
			& = \Gamma_{M, N} (f) \ab\big((v \otimes u)_n (a \otimes b)) \\
			& = (v \otimes u)_n \Gamma_{M, N} (f) (a \otimes b) \\
			& = (v \otimes u)_n \beta_{\Gamma} \ab\big(f(a \otimes b)) = (v \otimes u)_n \beta_{\Gamma}(m).
		\end{align*}
		Thus $\beta_\Gamma$ is indeed a $V\otimes U$-module homomorphism.
		Therefore $\Gamma = \mathcal{H}(\beta_{\Gamma})$.
		
		This completes the proof of surjectivity.
		
		\smallskip
		\noindent\textbf{Part 3: Equivalence for isomorphisms.}
		It is obvious that a $V \otimes U$-isomorphism $\beta \colon X \to Z$ induces a natural isomorphism
		$\Theta(\beta)$.
		Conversely, let $ \beta \colon X \to Z $ be a $V \otimes U$-homomorphism such that $\Theta(\beta)$
		is a natural isomorphism. We show that $\beta$ is itself a $V\otimes U$-module isomorphism.
		
		First, surjectivity: Take an arbitrary element $z \in Z$.
		Since $Z \in \mathcal{G}_{\mathcal{A}, \mathcal{B}}$,
		Lemma~\ref{lem:pure-tensor-presentation} guarantees the existence of an object $M \otimes N$ of $\mathcal{A} \odot \mathcal{B}$,
		an element $a \otimes b \in M \otimes N$, and a $V \otimes U$-module homomorphism $g \colon M \otimes N \to Z$ satisfying $g(a \otimes b) = z$.
		Since $\Theta(\beta)$ is a natural isomorphism, there exists a $V \otimes U$-module homomorphism $f \colon M \otimes N \to X$ with $\beta f = g$.
		Set $m = f(a \otimes b)$. Then
		\[
		\beta(m) = \beta(f(a \otimes b)) = g(a \otimes b) = z,
		\]
		which shows that $\beta$ is surjective.
		
		Now let $X_0 = \Ker(\beta)$. Corollary~\ref{cor:generated-category-closure} gives
		$X_0 \in \mathcal{G}_{\mathcal{A}, \mathcal{B}}$.
		For any $M \otimes N \in \mathcal{A}\odot\mathcal{B}$ and any $V \otimes U$-homomorphism
		$f \colon M \otimes N \to X_0$, we have $\beta f = 0$.
		Because $\Theta(\beta)$ is a natural isomorphism,
		this forces $f = 0$.
		Consequently, $X_0 = 0$.
		
		Therefore $\beta$ is a $V\otimes U$-module isomorphism.

	\end{proof}
	
	\vspace{0.2cm}
	
	\begin{proposition} \label{prop:generated-regular-values}
		Let $V$ be a vertex operator algebra, and let $W$ be a weak $V$-module.
		Let $\mathcal{A}$ and $\mathcal{B}$ be nonempty full subcategories of the category of weak $V$-modules such that both are closed under finite direct sums.
		Then
		\begin{align} \label{eq:generated-regular-values}
			\mathcal{R}_{\mathcal{A}, \mathcal{B}} (W)
			= \ab\{ \Y(a, y) b \; \bigm| \; M \in \mathcal{A}, \; N \in \mathcal{B}, \; a \in M, \; b \in N, \;\Y \in \mathcal{I}\binom{W}{M\; N} \}.
		\end{align}
	\end{proposition}
	
	\begin{proof}
		
		Let $\mathcal{Y} \in \mathcal{I}\binom{W}{M \; N} $ for some $M \in \mathcal{A}$, $N \in \mathcal{B}$.
		Then by Theorem~\ref{thm:hom-realization}, we have
		\[
		\mathcal{Y}(a, y) b = \psi^W_{M, N}(\mathcal{Y}) (a \otimes b) \in \psi^W_{M, N}(\mathcal{Y}) (M \otimes N) \subset \mathcal{R}_{\mathcal{A}, \mathcal{B}}(W),
		\]
		for any $a \in M$ and $b \in N$, where $\psi^W_{M, N}$ is the map defined in Theorem~\ref{thm:hom-realization}. Therefore, the right-hand side of \eqref{eq:generated-regular-values} is contained in the left-hand side.
		
		Take any $\phi(y) \in \mathcal{R}_{\mathcal{A}, \mathcal{B}}(W)$.
		By Lemma~\ref{lem:pure-tensor-presentation}, there exists a weak $V \otimes V$-module $M \otimes N \in \mathcal{A} \odot \mathcal{B}$, an element $a \otimes b \in M \otimes N$,
		and a $V \otimes V$-module homomorphism
		\[
		\mathcal{J} \colon M \otimes N \longrightarrow \mathcal{R}(W)
		\]
		such that $\mathcal{J}(a \otimes b) = \phi(y)$. Applying Theorem~\ref{thm:hom-realization} once more,
		there exists an intertwining operator $\mathcal{J}(\cdot, y)$ of type $\mathcal{I} \binom{W}{M\; N}$ such that
		\[
		\mathcal{J}(a, y) b = \mathcal{J}(a \otimes b) = \phi(y).
		\]
		As a consequence, the left-hand side of \eqref{eq:generated-regular-values} is contained in the right-hand side.
		
		This completes the proof.
		
	\end{proof}
	
	\begin{corollary} \label{cor:generated-regular-hom-realization}
		
		Let $V$ be a vertex operator algebra, and let $\mathcal{A}$ and $\mathcal{B}$ be nonempty full subcategories of the category of weak $V$-modules, both closed under finite direct sums and submodules.
		Let $W$ be a weak $V$-module.
		The natural isomorphism $\psi$ from Theorem~\ref{thm:hom-realization} restricts to a natural isomorphism of bifunctors
		\[
		\mathcal{I}^{W}(-, -) \quad \text{and} \;\; \mathcal{H}_{\mathcal{R}_{\mathcal{A}, \mathcal{B}}(W)}
		\]
		from $\mathcal{A}^{\opp}\times\mathcal{B}^{\opp}$ to $\mathbf{Vect}$.
		
	\end{corollary}
	
	\begin{proof}
		Combining Theorem~\ref{thm:hom-realization} with the isomorphism
		\[
		\Hom_{V \otimes V} \ab\big(M \otimes N, \, \mathcal{R}_{\mathcal{A}, \mathcal{B}}(W)) \cong
		\Hom_{V \otimes V} \ab\big(M \otimes N, \, \mathcal{R}(W))
		\]
		which holds for every $M \otimes N \in \mathcal{A} \odot \mathcal{B}$,
		we obtain the desired result immediately.
	\end{proof}
	
	The following theorem can be regarded as the Yoneda lemma for intertwining operator spaces.
	\begin{theorem}\label{thm:intertwiner-yoneda}
		Let $V$ be a vertex operator algebra, and let $\mathcal{A}$ and $\mathcal{B}$ be nonempty full subcategories of the category of weak $V$-modules, both closed under finite direct sums and submodules.
		Let $W_1$ and $W_2$ be weak $V$-modules.
		Then we have the following canonical isomorphism of vector spaces:
		\[
		\Nat_{\mathcal{A}^{\opp}\times\mathcal{B}^{\opp}} \ab\big(\mathcal{I}^{W_1}, \; \mathcal{I}^{W_2}) \cong
		\Hom_{V \otimes V} \ab\big( \mathcal{R}_{ \mathcal{A}, \mathcal{B} }(W_1), \; \mathcal{R}_{ \mathcal{A}, \mathcal{B} }(W_2) ),
		\]
		where $\mathcal{I}^{W_1}$ and  $\mathcal{I}^{W_2}$ are viewed, by restriction, as bifunctors on  $\mathcal{A}^{\opp} \times \mathcal{B}^{\opp}$.
		Moreover, under this isomorphism, natural isomorphisms correspond precisely to $V\otimes V$-module isomorphisms.
	\end{theorem}
	
	\begin{proof}
		
		From Corollary~\ref{cor:generated-regular-hom-realization} and Lemma~\ref{lem:external-tensor-density},
		\begin{align*}
			& \Nat_{\mathcal{A}^{\opp}\times\mathcal{B}^{\opp}} \ab\big(\mathcal{I}^{W_1}, \; \mathcal{I}^{W_2}) \\
			& \cong \Nat_{\mathcal{A}^{\opp}\times\mathcal{B}^{\opp}} \ab\Big( \mathcal{H}_{\mathcal{R}_{ \mathcal{A}, \mathcal{B} }(W_1)}, \;
			\mathcal{H}_{\mathcal{R}_{ \mathcal{A}, \mathcal{B} }(W_2)} ) \\
			& \cong \Hom_{V \otimes V} \ab\big( \mathcal{R}_{ \mathcal{A}, \mathcal{B} }(W_1), \; \mathcal{R}_{ \mathcal{A}, \mathcal{B} }(W_2) ).\\
		\end{align*}
		The correspondence between natural isomorphisms and $V\otimes V$-module isomorphisms also follows directly from Lemma~\ref{lem:external-tensor-density}. This completes the proof.
	\end{proof}
	
	\begin{corollary}
		
		Let $V$ be a vertex operator algebra, and let $W_1$ and $W_2$ be weak $V$-modules.
		Let $\mathcal{A}$ and $\mathcal{B}$ be nonempty full subcategories of the category of weak $V$-modules, both closed under finite direct sums and submodules.
		Suppose that $V \in \mathcal{A}$ and $W_1, W_2 \in \mathcal{B}$.
		Then the bifunctors $\mathcal{I}^{W_1}$ and $\mathcal{I}^{W_2}$ on $\mathcal{A}^{\opp}\times\mathcal{B}^{\opp}$
		are naturally isomorphic if and only if $W_1$ and $W_2$ are isomorphic as $V$-modules.
		
	\end{corollary}
	
	\begin{proof}
		
		It is clear that an isomorphism $W_1 \cong W_2$ of $V$-modules  gives rise to an isomorphism of bifunctors $\mathcal{I}^{W_1} \cong \mathcal{I}^{W_2} $.
		
		Suppose the bifunctors $\mathcal{I}^{W_1}$ and $\mathcal{I}^{W_2}$ are isomorphic. Then by Theorem~\ref{thm:intertwiner-yoneda}, there exist $V \otimes V$-module isomorphisms
		\[
		f \colon \mathcal{R}_{\mathcal{A}, \mathcal{B} }(W_1) \longrightarrow \mathcal{R}_{\mathcal{A}, \mathcal{B} }(W_2), \quad \text{and} \quad
		g \colon \mathcal{R}_{\mathcal{A}, \mathcal{B} }(W_2) \longrightarrow \mathcal{R}_{\mathcal{A}, \mathcal{B} }(W_1)
		\]
		such that
		\[
		gf = \id_{\mathcal{R}_{\mathcal{A}, \mathcal{B} }(W_1)}, \quad \text{and} \quad fg = \id_{ \mathcal{R}_{\mathcal{A}, \mathcal{B} }(W_2) }.
		\]
		Since $V \in \mathcal{A}$, $W_1 \in \mathcal{B}$, and $Y_{W_1} \in \mathcal{I} \binom{W_1}{V \; W_1}$,  Proposition~\ref{prop:generated-regular-values} implies that
		\[
		w_1 = Y_{W_1}(\mathbf{1}, y) w_1 \in \mathcal{R}(W_1), \quad \text{for all} \; w_1 \in W_1.
		\]
		Consequently, $W_1 \subset \mathcal{R}_{\mathcal{A}, \mathcal{B} }(W_1) $. Similarly, we have $W_2 \subset \mathcal{R}_{\mathcal{A}, \mathcal{B} }(W_2) $.
		
		On the other hand, by identity \eqref{eq:regular-translation}, for every $w_1 \in W_1$,
		\[
		\odv{y}f(w_1) = \omega^L_0f(w_1) = f(\omega^L_0 w_1) = f \ab\big(\odv{y} w_1) = 0.
		\]
		Hence $f(w_1)$  is a constant sequence, so $f(w_1) \in W_2$. Consequently, $f(W_1) \subset W_2$. Similarly, we have $g(W_2) \subset W_1$.
		Therefore, the restriction $f |_{W_1} \colon W_1 \to W_2$ is a $V$-module isomorphism. This completes the proof.
	\end{proof}
	
	\begin{remark}
		Let  $V$ and $U$ be vertex operator algebras.
		Let $\mathcal{A}$ (resp.\ $\mathcal{B}$) be a nonempty full subcategory of the category of weak $V$-modules (resp.\ weak $U$-modules).
		Suppose that $\mathcal{A}$ and $\mathcal{B}$ are closed under finite direct sums and submodules.
		Let $\mathcal{A}_0$ (resp.\ $\mathcal{B}_0$) be a full subcategory of $\mathcal{A}$ (resp.\ $\mathcal{B}$) such that every object of $\mathcal{A}$ (resp.\ $\mathcal{B}$) is isomorphic to some object of $\mathcal{A}_0$ (resp.\ $\mathcal{B}_0$).
		Then
		\[
		\mathcal{G}_{\mathcal{A}, \mathcal{B}} = \mathcal{G}_{\mathcal{A}_0, \mathcal{B}_0}, \quad \text{and} \quad
		\mathcal{S}_{\mathcal{A}, \mathcal{B}}(X) = \mathcal{S}_{\mathcal{A}_0, \mathcal{B}_0}(X)
		\]
		for any weak $V \otimes U$-module $X$.
		Moreover, the conclusions of Lemma~\ref{lem:pure-tensor-presentation}, Lemma~\ref{lem:common-presentation-refinement}, and Lemma~\ref{lem:external-tensor-density} remain valid if we replace $\mathcal{A}$ by $\mathcal{A}_0$ and $\mathcal{B}$ by $\mathcal{B}_0$. Consequently, all results in this subsection hold under this replacement.
		
		Note that under our definition, $\mathcal{A}_0$ and $\mathcal{B}_0$ need not be closed under finite direct sums and submodules.
		
	\end{remark}
	
	\subsection{ Coend Realization  } \label{Section-coend}
	
	Let $\mathcal{A}$ and $\mathcal{B}$ be categories, and let $Z$ be an object of $\mathcal{B}$.
	Let $F\colon\mathcal{A} \times\mathcal{A}^{\opp} \to \mathcal{B}$ be a bifunctor.
	A \textbf{ dinatural transformation} $\eta \colon F \ddot{ \Rightarrow } Z$ from the bifunctor $F$ to the object $Z$ is a family of morphisms
	\[
	\{\eta_{A}\colon F(A, A) \longrightarrow Z\}_{A \in \mathcal{A} }
	\]
	such that for any morphism $f\colon A \to B$ in $\mathcal{A}$, the following diagram commutes:
	\[
	\begin{tikzcd}
		{F(A, B)} \arrow[r, "{F(A, f)}"] \arrow[d, "{ F(f, B) }"'] & {F(A, A)} \arrow[d, "{\eta_A}"] \\
		{F(B, B)} \arrow[r, "{\eta_B}"'] & Z
	\end{tikzcd}
	\]
	
	The \textbf{ coend} of the bifunctor $F\colon\mathcal{A} \times\mathcal{A}^{\opp} \to \mathcal{B}$
	is an object $\Coend(F) \in \mathcal{B}$ equipped with a dinatural transformation
	$\rho \colon F \ddot{ \Rightarrow } \Coend(F)$ satisfying the following universal property: for any object $Z \in \mathcal{B}$,
	and any dinatural transformation $\eta \colon F \ddot{ \Rightarrow } Z$, there exists a unique morphism $g \colon \Coend(F) \to Z $
	such that
	\[
	g \circ \rho_A = \eta_A, \quad \text{for all} \; A \in \mathcal{A}.
	\]
	We usually denote $\Coend(F)$ by the following integral notation:
	\[
	\Coend(F) = \int^{A \in \mathcal A}F(A, A).
	\]

	Let $\mathcal{A}$ and $\mathcal{B}$ be categories, and let $F\colon\mathcal{A} \times\mathcal{A}^{\opp} \longrightarrow \mathcal{B} $ be a bifunctor.
	For each object $Z \in \mathcal{B}$, let $h_Z=\Hom_{\mathcal{B}}(-,Z)$. We obtain the bifunctor
	\[
	h_Z\circ F^{\opp} =  \Hom_{\mathcal{B}} \ab\big(F(-, -), \; Z) \colon \mathcal{A}^{\opp}\times\mathcal{A} \longrightarrow \mathbf{Vect}.
	\]

	\begin{lemma} \label{lem:dinatural-to-natural}
		Let $\mathcal{A}$ and $\mathcal{B}$ be categories, and let $F\colon\mathcal{A} \times\mathcal{A}^{\opp} \longrightarrow \mathcal{B} $ be a bifunctor.
		Let $Z \in \mathcal{B}$, and let $\eta \colon F \ddot{ \Rightarrow } Z$ be a dinatural transformation.
		For each object $ (A, B) \in \mathcal{A}^{\opp} \times\mathcal{A} $,  define a map
		\[
		\widetilde{\eta}_{A, B} \colon \Hom_{\mathcal{A}}(A, B)
		\longrightarrow \Hom_{\mathcal{B}} \ab\big(F(A, B), \; Z)
		\]
		by
		\begin{equation} \label{eq:dinatural-to-natural}
			\widetilde{\eta}_{A, B}(f) := \eta_B \circ F(f, B) = \eta_A \circ F(A, f),
		\end{equation}
		for all $f \colon A \to B$,
		where the equality follows from the dinaturality of $\eta$.
		Then the collection $\{\widetilde{\eta}_{A, B}\}$ assembles into a natural transformation
		\[
		\widetilde{\eta} \colon \Hom_{\mathcal{A}}(-, -) \Rightarrow h_Z\circ F^{\opp}.
		\]
	\end{lemma}
	
	\begin{proof}
		
		To prove the naturality of $\widetilde{\eta}$, take an arbitrary morphism
		\[
		(u, v) \colon (A_0, B_0) \longrightarrow (A, B)
		\]
		in $ \mathcal{A}^{\opp} \times\mathcal{A}$. This means precisely that
		\[
		u \colon A \to A_0, \quad \text{and} \quad v \colon B_0 \to B
		\]
		are morphisms in $\mathcal{A}$.
		Since
		\[
		(u, v) = (u, \id_B) \circ (\id_{A_0}, v),
		\]
		it is enough to verify naturality separately for the two morphisms
		\[
		(u, \id_B) \colon (A_0, B) \longrightarrow (A, B) \quad \text{and } \quad (\id_{A_0}, v) \colon (A_0, B_0) \longrightarrow (A_0, B).
		\]
		We first check naturality along the morphism $(u, \id_B)$.
		Equivalently, we verify that the following diagram commutes:
		\[
		\begin{tikzcd}
			{\Hom_{\mathcal{A}}(A_0, B)} \arrow[r, "{\widetilde{\eta}_{A_0,B}}"] \arrow[d, "{\Hom_{\mathcal{A}}(u, B)  }"']
			& {\Hom_{\mathcal{B}} \ab\big(F(A_0, B), \; Z)} \arrow[d, "{\Hom_{\mathcal{B}}(F(u, B), \;  Z)}"] \\
			{ \Hom_{\mathcal{A}}(A, B) } \arrow[r, "{\widetilde{\eta}_{A, B}}"'] & \Hom_{\mathcal{B}} \ab\big(F(A, B), \; Z)
		\end{tikzcd}
		\]
		By the definition of $\widetilde{\eta}$ and the functoriality of $F$,
		for any $f \in \Hom_{\mathcal{A}}(A_0, B)$, we compute:
		\begin{align*}
			& \ab\big(\widetilde{\eta}_{A, B} \circ \Hom_{\mathcal{A}}(u, B)) (f) \\
			& = \widetilde{\eta}_{A, B} (fu) = \eta_{B} \circ F(fu, B) = \eta_{B} \circ F(f, B) \circ F(u, B) \\
			& = \Hom_{\mathcal{B}} \ab\big(F(u, B), \; Z) \ab\big(\eta_{B} \circ F(f, B)) \\
			& = \Hom_{\mathcal{B}} \ab\big(F(u, B), \; Z) \ab\big(\widetilde{\eta}_{A_0, B} (f)) \\
			& = \ab\big(\Hom_{\mathcal{B}} \ab(F(u, B), \; Z) \circ \widetilde{\eta}_{A_0, B}) (f)
		\end{align*}
		This establishes commutativity of the diagram.
		
		Next we verify naturality for $ (\id_{A_0}, v)$. Again, we check that the diagram below commutes
		\[
		\begin{tikzcd}
			{\Hom_{\mathcal{A}}(A_0, B_0)} \arrow[r, "{\widetilde{\eta}_{A_0,B_0}}"] \arrow[d, "{\Hom_{\mathcal{A}}(A_0, v)  }"']
			& {\Hom_{\mathcal{B}} \ab\big(F(A_0, B_0), \; Z)} \arrow[d, "{\Hom_{\mathcal{B}}(F(A_0, v), \;  Z)}"] \\
			{ \Hom_{\mathcal{A}}(A_0, B) } \arrow[r, "{\widetilde{\eta}_{A_0, B}}"'] & \Hom_{\mathcal{B}} \ab\big(F(A_0, B), \; Z)
		\end{tikzcd}
		\]
		By the definition of $\widetilde{\eta}$ and the functoriality of $F$, for any $f \in \Hom_{\mathcal{A}}(A_0, B_0)$, we have
		\begin{align*}
			& \ab\big(\widetilde{\eta}_{A_0, B} \circ \Hom_{\mathcal{A}}(A_0, v)) (f) \\
			& = \widetilde{\eta}_{A_0, B} (vf) = \eta_{A_0} \circ F(A_0, vf) = \eta_{A_0} \circ F(A_0, f) \circ F(A_0, v) \\
			& = \Hom_{\mathcal{B}} \ab\big(F(A_0, v), \; Z) \ab\big(\eta_{A_0} \circ F(A_0, f)) \\
			& = \Hom_{\mathcal{B}} \ab\big(F(A_0, v), \; Z) \ab\big(\widetilde{\eta}_{A_0, B_0} (f)) \\
			& = \ab\big(\Hom_{\mathcal{B}}(F(A_0, v), \; Z) \circ \widetilde{\eta}_{A_0, B_0}) (f)
		\end{align*}
		This diagram commutes as well, which finishes the proof of naturality.
		
	\end{proof}
	
	\begin{lemma} \label{lem:natural-to-dinatural}
		Let $\mathcal{A}$ and $\mathcal{B}$ be categories, and let $F\colon\mathcal{A} \times\mathcal{A}^{\opp} \to \mathcal{B} $ be a bifunctor.
		Let $Z \in \mathcal{B}$, and suppose
		\[
		\gamma \colon \Hom_{\mathcal{A}}(-, -) \Rightarrow h_Z\circ F^{\opp}
		\]
		is a natural transformation. For each $A \in \mathcal{A}$,   define the map
		\begin{equation} \label{eq:diagonal-identity-evaluation}
			\mathring{\gamma}_A = \gamma_{A, A}(\id_A) \colon F(A, A) \longrightarrow Z.
		\end{equation}
		Then the family of morphisms  $\{\mathring{\gamma}_A\}$ forms a dinatural transformation
		$\mathring{\gamma} \colon F \ddot{ \Rightarrow } Z$.
		
	\end{lemma}
	
	\begin{proof}
		
		To show that the family of morphisms  $\{\mathring{\gamma}_A\}$ forms a dinatural transformation,
		we need to check that for every morphism  $f \colon A \to B$ in $\mathcal{A}$, the equality
		\[
		\mathring{\gamma}_B \circ F(f, B) = \mathring{\gamma}_A \circ F(A, f)
		\]
		holds.  This follows from the naturality of $\gamma$.
		Indeed, the naturality of $\gamma$ provides the two commutative diagrams below:
		\[
		\begin{tikzcd}
			{\Hom_{\mathcal{A}}(B, B)} \arrow[r, "{\gamma_{B, B}}"] \arrow[d, "{\Hom_{\mathcal{A}}(f, B)  }"']
			& {\Hom_{\mathcal{B}} \ab\big(F(B, B), \; Z) } \arrow[d, "{\Hom_{\mathcal{B}}(F(f, B), \;  Z)}"] \\
			{ \Hom_{\mathcal{A}}(A, B) } \arrow[r, "{ \gamma_{A, B} }"'] & \Hom_{\mathcal{B}} \ab\big(F(A, B), \;  Z)
		\end{tikzcd}
		\]
		and
		\[
		\begin{tikzcd}
			{\Hom_{\mathcal{A}}(A, A)} \arrow[r, "{\gamma_{A, A}}"] \arrow[d, "{\Hom_{\mathcal{A}}(A, f)  }"']
			& {\Hom_{\mathcal{B}} \ab\big(F(A, A), \;  Z)} \arrow[d, "{\Hom_{\mathcal{B}} (F(A, f), \;  Z)}"] \\
			{ \Hom_{ \mathcal{A} }(A, B) } \arrow[r, "{ \gamma_{A, B} }"'] & \Hom_{\mathcal{B}} \ab\big(F(A, B), \;  Z)
		\end{tikzcd}
		\]
		From the first diagram,
		\begin{equation}\label{eq:diagonal-naturality}
			\begin{aligned}
				\mathring{\gamma}_B \circ F(f, B) & = \gamma_{B, B}(\id_B) \circ F(f, B) \\
				& = \Hom_{\mathcal{B}} \ab\big(F(f, B), \; Z) \ab\big(\gamma_{B, B} (\id_B)) \\
				& = \ab\big(\Hom_{\mathcal{B}}(F(f, B), \; Z) \circ \gamma_{B, B}) (\id_B) \\
				& = \ab\big(\gamma_{A, B} \circ \Hom_{\mathcal{A}}(f, B)) (\id_B) = \gamma_{A, B}(f).
			\end{aligned}
		\end{equation}
		Similarly, the second diagram yields
		\begin{align*}
			\mathring{\gamma}_A \circ F(A, f) & = \gamma_{A, A}(\id_A) \circ F(A, f) \\
			& = \Hom_{\mathcal{B}} \ab\big(F(A, f), \; Z) \ab\big(\gamma_{A, A} (\id_A)) \\
			& = \ab\big(\Hom_{\mathcal{B}}(F(A, f), \; Z) \circ \gamma_{A, A}) (\id_A) \\
			& = \ab\big(\gamma_{A, B} \circ \Hom_{\mathcal{A}}(A, f)) (\id_A) = \gamma_{A, B}(f).
		\end{align*}
		Therefore,
		\[
		\mathring{\gamma}_B \circ F(f, B) = \mathring{\gamma}_A \circ F(A, f),
		\]
		which completes the verification.
		
	\end{proof}
	
	\begin{corollary} \label{cor:natural-dinatural-bijection}
		Let $\mathcal{A}$ and $\mathcal{B}$ be categories, let $Z \in \mathcal{B}$, and let $F\colon\mathcal{A} \times\mathcal{A}^{\opp} \to \mathcal{B} $ be a bifunctor.
		Define
		\[
		\Dinat(F, \; Z) = \{ \text{ dinatural transformations from $F$ to $Z$ } \}.
		\]
		Then the assignment
		\[
		\Nat \ab\big(\Hom_{\mathcal{A}}(-, -), \; h_Z\circ F^{\opp})
		\longrightarrow \Dinat(F, \; Z), \quad \gamma \longrightarrow \mathring{\gamma}
		\]
		is bijective.
	\end{corollary}
	
	\begin{proof}
		Let $\gamma$ and $\theta$ be natural transformations such that $\mathring{\gamma} = \mathring{\theta}$.
		Then for any morphism in $f\colon A \to B$ in $\mathcal{A}$, by identity \eqref{eq:diagonal-naturality}, we have
		\[
		\gamma(f) = \mathring{\gamma}_B\circ F(f, B) = \mathring{\theta}_B \circ F(f, B) = \theta(f).
		\]
		Hence this assignment is injective.
		
		To establish surjectivity of this assignment,
		take any $\eta \in \Dinat(F, \; Z).$
		Using \eqref{eq:dinatural-to-natural} and \eqref{eq:diagonal-identity-evaluation},  we obtain
		\[
		\ab\big(\mathring{\widetilde{\eta}})_A = \widetilde{\eta}_{A, A}(\id_A) = \eta_A \circ F(A, \id_A) = \eta_A
		\]
		for any $A \in \mathcal{A}$.
		This completes the proof.
	\end{proof}
	
	\begin{lemma} \label{lem:vertex-operator-intertwiners}
		Let $V$ be a vertex operator algebra, and let $(M, Y_M)$ and $(N, Y_N)$ be weak $V$-modules.
		Define the linear map
		\[
		\tau_{M, N}\colon \Hom_V(M, N) \longrightarrow \mathcal{I}\binom{N}{V\,M}, \qquad f \mapsto [\id_V, \id_M, f] \circ Y_M.
		\]
		Then $\tau_{M, N}$ is a linear isomorphism.
	\end{lemma}
	
	\begin{proof}
		This Lemma is essentially due to \cite[Remark 4.11]{Li02b}. We reproduce a proof for completeness.
		We first observe that $\tau_{M, N}$ is clearly injective and linear.
		Now let $\mathcal{Y} \in \mathcal{I}\binom{N}{V\,M}$. Recall that $L(-1) \mathbf{1} = 0$. Hence, for all $a \in M$,
		\[
		\odv{y}\mathcal{Y}(\mathbf{1}, y) a = \mathcal{Y}(L(-1) \mathbf{1}, y) a = 0.
		\]
		Thus $\mathcal{Y}(\mathbf{1}, y) a$ is independent of $y$.
		We may therefore define a linear map $f\colon M \to N$ by
		\[
		f(a) = \mathcal{Y}(\mathbf{1}, y) a, \quad a \in M.
		\]
		We claim that $f$ is a $V$-module homomorphism and satisfies
		$[\id_V, \id_M, f] \circ Y_M = \mathcal{Y}$.
		
		Take arbitrary $u \in V$ and $a \in M$. Since $Y(u, x) \mathbf{1} \in V\dbrack{x}$, the weak commutativity of the intertwining operator $\mathcal{Y}$ yields
		\[
		Y_N(u, x) f(a) = Y_N(u, x) \mathcal{Y}(\mathbf{1}, y) a
		= \mathcal{Y}(\mathbf{1}, y) Y_M(u, x) a
		= f \ab\big(Y_M(u, x) a).
		\]
		This verifies that $f$ is a $V$-module homomorphism.
		
		Next, choose $q \in \N$ such that $x^q Y_M(u, x) a \in M\dbrack{x}$. By the weak associativity of $\mathcal{Y}$, we have
		\[
		(z+y)^q Y_N(u, z+y) \mathcal{Y}(\mathbf{1}, y) a
		= (z+y)^q \mathcal{Y} \ab\big(Y_V(u, z) \mathbf{1}, \,y) a.
		\]
		Setting $z = 0$ and cancelling the common factor $y^q$ gives
		\[
		\mathcal{Y}(u, y) a = Y_N(u, y) f(a) = f \ab\big(Y_M(u, y) a).
		\]
		Consequently, $[\id_V, \id_M, f] \circ Y_M = \mathcal{Y}$, which proves the surjectivity of $\tau_{M, N}$.
		
		Therefore $\tau_{M, N}$ is a linear isomorphism.  This completes the proof.
	\end{proof}
	
	Let $V$ be a vertex operator algebra.
	Let $\mathcal{A}$ be a nonempty full subcategory of the category of restricted $V$-modules that is closed under finite direct sums, submodules, and restricted duals. Consequently, $\mathcal{A}$ is automatically closed under quotients,  and
	\[
	\mathcal{A} \odot \mathcal{A} = \ab\{ M \otimes N^\prime \; \big | \; M, N \in \mathcal{A} \}.
	\]
	
	\begin{lemma} \label{lem:hom-kernel-realization}
		Let $V$ and $\mathcal{A}$ be as above. Let $(M, Y_M)$ and $(N, Y_N)$ be objects of $\mathcal{A}$.
		Define a linear map
		\[
		\varphi_{M, N} \colon \Hom_V(M, N) \longrightarrow \Hom_{V \otimes V} \ab\big(M \otimes N^\prime, \; \mathcal{R}_{\mathcal{A}, \mathcal{A} }(V^\prime))
		\]
		as the composition of the following sequence of maps:
		\[
		\Hom_V(M, N) \xrightarrow{\tau_{M, N}} \mathcal{I}\tbinom{N}{V\;M} \xrightarrow{\Omega_0} \mathcal{I}\tbinom{N}{M\;V}
		\xrightarrow{T_0} \mathcal{I}\tbinom{V^\prime}{M\;N^\prime}
		\xrightarrow{\psi^{V^\prime}_{M, N^\prime}} \Hom_{V \otimes V} \ab\big(M \otimes N^\prime, \; \mathcal{R}_{\mathcal{A}, \mathcal{A} }(V^\prime)),
		\]
		where, $\tau_{M, N}$ is the map defined in Lemma~\ref{lem:vertex-operator-intertwiners}, and $\psi^{V^\prime}_{M, N^\prime}$ is the map defined in Theorem~\ref{thm:hom-realization}.
		Then $\varphi_{M, N}$ is a linear isomorphism.
		Moreover, these isomorphisms assemble into an isomorphism of natural transformations
		\[
		\varphi \colon \Hom_{V}(-, -) \Rightarrow \mathcal{H}_{\mathcal{R}_{\mathcal{A}, \mathcal{A} }(V^\prime)}(-, -^\prime)=\Hom_{V \otimes V}\big(- \otimes-^\prime, \; \mathcal{R}_{\mathcal{A}, \mathcal{A} }(V^\prime) \big) ,
		\]
		where both sides are regarded as bifunctors from $\mathcal{A}^{\opp} \times \mathcal{A}$ to $\mathbf{Vect}$.
	\end{lemma}
	
	\begin{proof}
		
		Since $\varphi_{M, N}$ is a composition of isomorphisms, it is therefore a linear isomorphism.
		
		The naturality of $\varphi$  is equivalent to the following condition:   for all $V$-module homomorphisms $\alpha\colon M_0 \to M$ and $\beta\colon N \to N_0$,
		the outer rectangle in the diagram below commutes.This follows by pasting together the four adjacent commuting rectangular subdiagrams.
		\[
		\begin{tikzcd}[font=\fontsize{10pt}{11pt}\selectfont, column sep=small]
			\Hom_{V}(M, N) \arrow[r, "{\tau_{M, N}}"] & \mathcal{I}\binom{N}{V\; M} \arrow[r, "\Omega_0"] & \mathcal{I}\binom{N}{M\; V} \arrow[r, "T_0"] & \mathcal{I}\binom{V^\prime}{M\; N^\prime} \arrow[r, "{\psi^{V^\prime}_{M,N^\prime}}"] & \Hom_{V \otimes V} \ab\big(M \otimes N^\prime, \; \mathcal{R}_{\mathcal{A}, \mathcal{A}}(V^\prime)) \\
			\Hom_{V}(M_0, N_0) \arrow[r, "{\tau_{M_0, N_0}}"] & \mathcal{I}\binom{N_0}{V\; M_0} \arrow[r, "\Omega_0"] & \mathcal{I}\binom{N_0}{M_0\; V}  \arrow[r, "T_0"] & \mathcal{I}\binom{V^\prime}{M_0\; N^\prime_0}  \arrow[r, "{\psi^{V^\prime}_{M_0, N_0^\prime}}"] & \Hom_{V \otimes V} \ab\big(M_0 \otimes N_0^\prime, \; \mathcal{R}_{\mathcal{A}, \mathcal{A}}(V^\prime))
			\arrow[from=1-1, to=2-1, "{\Hom_V(\alpha, \beta)}"]
			\arrow[from=1-2, to=2-2, "{[\id, \alpha, \beta] \circ-}"]
			\arrow[from=1-3, to=2-3, "{[\alpha, \id, \beta] \circ-}"]
			\arrow[from=1-4, to=2-4, "{[\alpha, \beta^\prime, \id] \circ-}"]
			\arrow[from=1-5, to=2-5, "{ \Hom_{V \otimes V}(\alpha \otimes \beta^\prime,\; \mathcal{R}_{\mathcal{A}, \mathcal{A}}(V^\prime)) }"]
		\end{tikzcd}
		\]
		Commutativity of the second rectangular subdiagram follows from Identity~\eqref{eq:skew-symmetry-naturality}, that of the third from Identity~\eqref{eq:contragredient-naturality},
		and that of the fourth from the naturality of $\psi^{V^\prime}$.
		
		Given any $V$-module homomorphism $f\colon M \to N$, any $v \in V$ and $a_0 \in M_0$,
		by the definition of $\tau$ and the fact that $\alpha$ and $\beta$ are $V$-module homomorphisms, we have
		\begin{align*}
			\ab\big([\id, \alpha, \beta] \circ \ab\big(\tau_{M, N} (f))) (v, y) a_0
			& = \ab\big([\id, \alpha, \beta] \circ \ab\big([\id, \id, f] \circ Y_M)) (v, y) a_0\\
			& = \beta \ab\big(\ab\big([\id, \id, f] \circ Y_M) (v, y) \alpha(a_0)) \\
			& = \beta f \ab\big(Y_M(v, y) \alpha(a_0)) = \beta f \alpha \ab\big(Y_{M_0}(v, y) a_0) \\
			& = \ab\big(\tau_{M_0, N_0}(\beta f \alpha)) (v, y) a_0 \\
			& = \ab\big(\tau_{M_0, N_0} \ab\big(\Hom_V(\alpha, \beta) (f))) (v, y) a_0.
		\end{align*}
		It follows that
		\[
		[\id, \alpha, \beta] \circ \tau_{M, N} = \tau_{M_0, N_0} \circ \Hom_V(\alpha, \beta),
		\]
		which establishes the commutativity of the first rectangular subdiagram and completes the proof.
		
	\end{proof}

	\begin{theorem} \label{thm:regular-coend-realization}
		Let $V$ and $\mathcal{A}$ be as above. Let $\mathcal{B} $ be a full subcategory of the category of weak $V \otimes V$-modules that contains
		$\mathcal{G}_{\mathcal{A}, \mathcal{A}}$.
		Consider the bifunctor
		\[
		F = - \otimes -^\prime \colon \mathcal{A} \times \mathcal{A}^{\opp} \longrightarrow \mathcal{B},
		\quad (M, N) \longrightarrow M \otimes N^\prime.
		\]
		Then the coend of $F$ exists and is given by
		\[
		\int^{M} M \otimes M^\prime = \mathcal{R}_{\mathcal{A}, \mathcal{A}}(V^\prime).
		\]
	\end{theorem}
	
	\begin{proof}
		
		Let $\varphi$ be the natural transformations defined in Lemma~\ref{lem:hom-kernel-realization}.
		For each $M \in \mathcal{A}$, define a  morphism
		\[
		\pi_M = \varphi_{M, M}(\id_M) \colon M\otimes M^\prime \to \mathcal{R}_{\mathcal{A}, \mathcal{A}}(V^\prime).
		\]
		By Lemma~\ref{lem:natural-to-dinatural},
		the family $\{\pi_M\}_{M \in \mathcal{A} }$ forms a dinatural transformation $\pi \colon F \ddot{ \Rightarrow } \mathcal{R}_{\mathcal{A}, \mathcal{A}}(V^\prime)$.
		We next show that this dinatural transformation satisfies the universal property of the coend of $F$.
		
		Let $Z \in \mathcal{B}$,  and let $\eta: F \ddot{ \Rightarrow } Z$ be any dinatural transformation.
		Define, for each $(M, N) \in \mathcal{A}^{\opp} \times \mathcal{A}$, a map
		\[
		\gamma_{M, N}\colon \Hom_{\mathcal{A}}(M, N) \longrightarrow \Hom_{\mathcal{B}} \ab\big(M \otimes N^\prime, \; Z)
		\]
		by
		\[
		\gamma_{M, N}(f) := \eta_N \circ (f \otimes \id_{N^\prime}).
		\]
		Then by Lemma~\ref{lem:dinatural-to-natural}, these maps assemble into a natural transformation
		\[
		\gamma \colon \Hom_{\mathcal{A}}(-, -) \Rightarrow \mathcal{H}_Z(-, -^\prime).
		\]
		Consequently, by composing with $\varphi^{-1}$,  we obtain a natural transformation
		\[
		\gamma \circ \varphi^{-1} \colon \mathcal{H}_{\mathcal{R}_{\mathcal{A}, \mathcal{A}}(V^\prime)}(-, -^\prime) \Rightarrow \mathcal{H}_Z(-, -^\prime).
		\]
		Since $M \cong (M^\prime)^\prime$ for any $M \in \mathcal{A}$, it follows from Lemma~\ref{lem:external-tensor-density} that there exists a unique  $V \otimes V$-homorphism $f \colon \mathcal{R}_{\mathcal{A}, \mathcal{A}}(V^\prime) \to Z $   such that the following diagram commutes for all
		$M, N \in \mathcal{A}$:
		\[
		\begin{tikzcd}
			\Hom_{\mathcal{A}}(M, N)
			\arrow[r, "{\varphi_{M, N}}"]
			\arrow[dr, "{\gamma_{M, N}}"']
			& \Hom_{\mathcal{B}}(M \otimes N^\prime, \; \mathcal{R}_{\mathcal{A}, \mathcal{A}}(V^\prime))
			\arrow[d, "{\Hom_{\mathcal{B}}(M \otimes N^\prime,  \; f)}"]
			\\
			& \Hom_{\mathcal{B}}(M\otimes N^\prime, \; Z)
		\end{tikzcd}
		\]
		Taking $N = M$ gives the commutative diagram
		\[
		\begin{tikzcd}[column sep=large, row sep=normal]
			M \otimes M^\prime
			\arrow[r, "{\varphi_{M, M}(\id_M)}"]
			\arrow[dr, "{\gamma_{M,M}(\id_M) }"']
			& \mathcal{R}_{\mathcal{A}, \mathcal{A}}(V^\prime)
			\arrow[d, "f"]
			\\
			& Z
		\end{tikzcd}
		\]
		Since $\pi_M = \varphi_{M, M}(\id_M) $ and $ \gamma_{M, M}(\id_M) = \eta_M $, we obtain
		\begin{align} \label{eq:coend-factorization}
			f \circ \pi_M = \eta_M \quad \text{for all} \; M \in \mathcal{A}.
		\end{align}
		Finally, because $ \mathcal{R}_{\mathcal{A}, \mathcal{A}}(V^\prime) $ is generated by $ \mathcal{A} \odot \mathcal{A} $,
		the morphism $f$ satisfying \eqref{eq:coend-factorization} is necessarily unique. This verifies the universal property of the coend, and thus the proof is complete.
		
	\end{proof}
	
	For $\mathcal{A} = \mathfrak{M}^r_V$, Proposition~\ref{prop:generated-category-characterizations}(3) gives
	\[
	\mathcal{R}_{\lr}(V^\prime)
	= \mathcal{R}_{\mathfrak{M}^r_V, \mathfrak{M}^r_V}(V^\prime).
	\]
	
	\begin{corollary}\label{cor:restricted-coend-realization}
		The coend over all restricted
		$V$-modules exists in $\mathfrak{M}^w_{V\otimes V}$ and satisfies
		\[
		\mathcal{R}_{\lr}(V')
		\cong\int^{M\in\mathfrak{M}^r_V}M\otimes M'.
		\]
	\end{corollary}
	
	\begin{proof}
		Apply Theorem~\ref{thm:regular-coend-realization} with
		$\mathcal{A}=\mathfrak{M}^r_V$ and use the preceding identification.
	\end{proof}

	\begin{remark}
		Theorem~\ref{thm:regular-coend-realization} extends Gui--Zhang's VOA realization of Lyubashenko's coend formula (see \cite[Corollary~2.34]{GZ25} and \cite[Theorem~0.1]{GZ26b}) beyond the $C_2$-cofinite setting. Their proofs are built in the context of the Deligne tensor product of finite abelian categories. By comparison, our proof uses a Yoneda-type lemma and does not assume such finiteness.
	\end{remark}
	
	\section{Relation to Li's Regular Representations} \label{sec:li-comparison}
	
	Throughout this section, let $\xi$ be a fixed nonzero complex number. We adopt the notation
	\[
	\Log \xi = \ln | \xi | + \operatorname{i} \Arg \xi, \quad \text{where } 0 \leq \Arg \xi < 2\pi.
	\]
	For $\alpha \in \C $, we then define
	\[
	{\xi}^{\alpha} = e^{\alpha \Log \xi}.
	\]
	\begin{definition} \cite{Li02}
		Let $V$ be a vertex operator algebra, and let $W$ be a weak $V$-module.
		We denote by $\mathfrak{D}_{\xi}(W)$ the set consisting of all linear functionals $f \in W^*$ satisfying the following condition:
		for every $u \in V$, there exist $p, q \in \N$ (depending on $f$ and $u$)  such that
		\begin{equation}\label{eq:li-functional-condition}
			x^p(x-\xi)^qY_{W^*}(u, x) f \in W^*\dbrack{x}.
		\end{equation}
		Equivalently, the condition is that
		\[
		x^p(x-\xi)^q \ab\big \langle f, \; Y^\circ_W(u, x) w\rangle \in \C[x]
		\]
		for all $w \in W$.
		
	\end{definition}
	
	\begin{definition}
		Let $V$ be a vertex operator algebra, and let $W$ be a weak $V$-module.
		Suppose $u \in V$ and $f \in \mathfrak{D}_{\xi}(W)$.   Choose  $p, q \in \N$ such that  \eqref{eq:li-functional-condition}  holds.
		Define
		\[
		Y^{R}_{ \mathfrak{D}_{\xi}(W) }(u, x) f, \quad
		Y^{L}_{ \mathfrak{D}_{\xi}(W) }(u, x) f \in W^*\dparen{x}
		\]
		by the following formulas
		\[
		\ab\big \langle Y^{R}_{ \mathfrak{D}_{\xi}(W) }(u, x) f, \; w\rangle
		= (-\xi+x)^{-q} \ab\big((x-\xi)^q \ab\big \langle f, \; Y^\circ_W(u, x) w\rangle),
		\]
		and
		\[
		\ab\big \langle Y^{L}_{ \mathfrak{D}_{\xi}(W) }(u, x) f, \; w\rangle
		= (\xi+x)^{-p} \ab\big((x+\xi)^p \ab\big \langle f, \; Y^\circ_W(u, x+\xi) w\rangle)
		\]
		respectively,  for all $w \in W$.
	\end{definition}
	
	We point out that the definitions of $Y^{R}_{ \mathfrak{D}_{\xi}(W) }$ and $Y^{L}_{ \mathfrak{D}_{\xi}(W) }$
	are independent of the choice of $p, q \in \N$ satisfying \eqref{eq:li-functional-condition}.
	
	\begin{proposition}\cite{Li02} \label{prop:li-commuting-actions}
		Let $V$ be a vertex operator algebra and $W$ a weak $V$-module. Then:
		\begin{enumerate}[label=(\arabic*),ref=\arabic*]
			\item Both pairs
			$\ab\big(\mathfrak{D}_{\xi}(W), Y^{L}_{ \mathfrak{D}_{\xi}(W) })$ and $\ab\big(\mathfrak{D}_{\xi}(W), Y^{R}_{ \mathfrak{D}_{\xi}(W) })$
			are weak $V$-modules.
			
			\item The actions of $Y^{R}_{\mathfrak{D}_{\xi}(W)}(-, x)$ and $Y^{L}_{\mathfrak{D}_{\xi}(W)}(-, x)$ on $\mathfrak{D}_{\xi}(W)$ commute.
			Consequently, $\mathfrak{D}_{\xi}(W)$ is a weak $V \otimes V$-module with the vertex operator map $Y_{\mathfrak{D}_{\xi}(W)}$ defined
			by
			\[
			Y_{\mathfrak{D}_{\xi}(W)}(u \otimes v, x) = Y^{R}_{\mathfrak{D}_{\xi}(W)}(u, x) Y^{L}_{\mathfrak{D}_{\xi}(W)}(v, x) \quad \text{for all} \; u, v \in V.
			\]
		\end{enumerate}
		
	\end{proposition}
	
	Let $W = \bigoplus_{\lambda \in \C} W_{(\lambda)}$ be a generalized $V$-module, and set
	\[
	\overline{W} = \prod_{\lambda \in \C} W_{(\lambda)}.
	\]
	For each $\lambda \in \C$, we let $\pi_\lambda$ denote the canonical projection $\overline{W} \to W_{(\lambda)}$.
	For $u \in V$, let $Y_{\overline{W}}(u, x)$ denote the componentwise action of $Y_W(u, x)$ on $\overline{W}$.

	\begin{definition}
		Let $V$ be a vertex operator algebra, and let $M, N, W$ be generalized $V$-modules.
		A \textbf{ $P(\xi)$-intertwining map } of type $\binom{W}{M \; N}$ is a linear map
		\[
		I \colon M \otimes N \longrightarrow \overline{W}
		\]
		satisfying the following conditions:
		\begin{enumerate}[label=(\arabic*),ref=\arabic*]
			\item \textbf{Lower truncation}: for any $a \in M, b \in N$ and $\lambda \in \C$,
			\[
			\pi_{\lambda-n} \ab\big(I(a \otimes b)) = 0 \quad \text{for} \; n \in \N \; \text{sufficiently large.}
			\]
			
			\item \textbf{The Jacobi identity}: for all $a \in M, b \in N$ and $u \in V$,
			\begin{align*} \label{eq:p-intertwining-jacobi}
				& x_0^{-1} \delta \ab(\frac{x_1- \xi}{x_0}) Y_{\overline{W}}(u, x_1) I(a \otimes b)
				-x_0^{-1} \delta \ab(\frac{-\xi+x_1}{x_0}) I \ab\big(a \otimes Y_{N}(u, x_1) b) \\
				& = {\xi}^{-1} \delta \ab(\frac{x_1-x_0}{\xi}) I \ab\big(Y_{M}(u, x_0) a \otimes b).
			\end{align*}
		\end{enumerate}
		Let $ \mathcal{I}_{P(\xi)}\binom{W}{M \; N} $ denote the space of $P(\xi)$-intertwining maps of type $\binom{W}{M \; N}$.
		
	\end{definition}
	
	Let $M, N, W$ be generalized $V$-modules, and suppose $\Y \in \mathcal{I}\binom{W}{M \; N}$.
	For $a \in M_{(\lambda)}, b \in N_{(\mu)}$ with $\lambda, \mu \in \C$, and $n \in \C$, we have
	\[
	a_{n; k}b \in W_{(\lambda +\mu -n-1)}
	\]
	for all $k \in \N$ \cite{HLZ10II, YZ26a}. It follows that
	\[
	\Y_{\xi}(a \otimes b) := \ab\big(\Y(a, y) b)_{ \big | y^{-n-1} (\log y)^k = {\xi}^{-n-1} (\Log \xi)^k }
	\]
	is a well-defined element in $\overline{W}$. Consequently, $\Y_{\xi} \in \mathcal{I}_{P(\xi)}\binom{W}{M \; N}.$
	The property below shows that every $P(\xi)$-intertwining map is obtained in this fashion.
	
	\begin{proposition} \cite{HLZ10III} \label{prop:intertwiner-evaluation}
		Let $V$ be a vertex operator algebra, and let $M, N, W$ be generalized $V$-modules.
		Then the map
		\[
		\mathcal{I}\binom{W}{M \; N} \longrightarrow \mathcal{I}_{P(\xi)}\binom{W}{M \; N}, \quad \Y \mapsto \Y_{\xi}
		\]
		is a linear isomorphism.
		
	\end{proposition}
	
	Let $M, N$ be generalized $V$-modules, and let $W$ be a restricted $V$-module.
	Then $\overline{W^\prime}$ is canonically isomorphic to $W^*$.
	Suppose that $I \in \mathcal{I}_{P(\xi)}\binom{W^\prime}{M \; N}$.
	Let $a \in M$ and $b \in N$.
	For each $u \in V$, choose $k \in \N$ such that $x^kY_{M}(u, x) a \in M\dbrack{x}$. Then, by the Jacobi identity for $I$, we have
	\[
	(x- \xi)^k Y_{W^*}(u, x) I(a \otimes b) = (-\xi + x)^kI \ab\big(a \otimes Y_{N}(u, x) b) \in W^*\dparen{x}.
	\]
	Consequently, $I(a \otimes b) \in \mathfrak{D}_{\xi}(W)$ for all $a \in M$ and $b \in N$.
	Moreover, the following proposition holds. It extends the ordinary module case considered in \cite[thoerem 4.5]{Li02} to the restricted module case, and its proof is exactly the same.
	
	\begin{proposition} \label{prop:p-intertwining-li-hom}
		Let $V$ be a vertex operator algebra, and let $W$ be a restricted $V$-module.
		For any generalized $V$-modules $M$ and $N$, the linear map
		\[
		\widetilde{\gamma}_{M, N}\colon \mathcal{I}_{P(\xi)} \binom{W^\prime}{M \; N} \longrightarrow \Hom_{V \otimes V} \ab\big(M \otimes N, \; \mathfrak{D}_{\xi}(W))
		\]
		defined by
		\[
		\widetilde{\gamma}_{M, N}(I) (a \otimes b) = I(a \otimes b) \quad \ab\big(I \in \mathcal{I}_{P(\xi)} \binom{W^\prime}{M \; N}, a \in M, b \in N),
		\]
		is a linear isomorphism.
	\end{proposition}
	
	\begin{corollary} \label{cor:li-intertwiner-hom}
		Let $V$ be a vertex operatro algebra, and let $W$ be a restricted $V$-module.
		For any generalized $V$-modules $M$ and $N$, the linear map
		\begin{align*}\label{eq:li-intertwiner-hom-map}
			\gamma_{M, N}\colon \mathcal{I} \binom{W^\prime}{M \; N} \longrightarrow \Hom_{V \otimes V} \ab\big(M \otimes N, \; \mathfrak{D}_{\xi}(W))
		\end{align*}
		defined by
		\[
		\gamma_{M, N}(\Y) (a \otimes b) = \Y_{\xi}(a \otimes b) \quad \ab\big(\Y \in \mathcal{I}\binom{W^\prime}{M \; N}, a \in M, b \in N),
		\]
		is a linear isomorphism.
		Moreover, these isomorphisms assemble into a natural isomorphism
		\[
		\gamma \colon \mathcal{I}^{W^\prime}(-, -) \Rightarrow
		\mathcal{H}_{\mathfrak{D}_{\xi}(W)},
		\]
		where  both sides are regarded as bifunctors from $(\mathfrak{M}^g_V)^{\opp}\times(\mathfrak{M}^g_V)^{\opp}$ to $\mathbf{Vect}$.
	\end{corollary}
	
	\begin{proof}
		It follows directly from Proposition~\ref{prop:intertwiner-evaluation}  and~\ref{prop:p-intertwining-li-hom} that
		$\gamma_{M, N}$ is a linear isomorphism. Naturality is verified in a way similar to that in Theorem~\ref{thm:hom-realization}.
	\end{proof}
	
	\begin{theorem} \label{thm:li-regular-comparison}
		Let $V$ be a vertex operator algebra,
		and let $W$ be a restricted $V$-module. Let $\mathcal{A} = \mathfrak{M}^g_V$ be the category of generalized $V$-modules.
		Then the linear map
		\[
		\ev_{\xi} \colon \mathcal{R}_{\A, \A} (W^\prime) \longrightarrow \mathcal{S}_{\A, \A} \ab\big(\mathfrak{D}_{\xi}(W))
		\]
		defined by
		\[
		\ev_{\xi} \ab\big(\Y(a, y) b) = \Y_{\xi}(a \otimes b)
		\]
		is a linear isomorphism, where $\Y$ is an intertwining operator of type $\binom{W^\prime}{M \; N}$
		for some generalized $V$-modules $M$ and $N$, with $a \in M$ and $b \in N$.
		In particular, the unique largest bi-generalized $V\otimes V$-submodule of $\mathcal{R}(W')$ is isomorphic to that of $\mathfrak{D}_{\xi}(W)$.
		
	\end{theorem}
	
	\begin{proof}
		By Corollary~\ref{cor:generated-regular-hom-realization} and \ref{cor:li-intertwiner-hom}, we obtain natural isomorphisms of bifunctors
		$ \mathcal{A}^{\opp} \times \mathcal{A}^{\opp}$  to $\mathbf{Vect}$:
		\[
		\mathcal{H}_{\mathcal{R}_{\A, \A}(W^\prime)} \cong
		\mathcal{H}_{\mathcal{S}_{\A, \A} (\mathfrak{D}_{\xi}(W))}.
		\]
		By Lemma~\ref{lem:external-tensor-density}, this natural isomorphism of bifunctors gives a $V \otimes V$-module isomorphism
		\[
		\ev_{\xi} \colon \mathcal{R}_{\A, \A} (W^\prime) \longrightarrow \mathcal{S}_{\A, \A} \ab\big(\mathfrak{D}_{\xi}(W))
		\]
		such that, for any $M, N \in \mathcal{A}$, the diagram below is commutative:
		\[
		\begin{tikzcd}
			\mathcal{I}\binom{W^\prime}{M \; \;N}
			\arrow[r, "{\psi^{W^\prime}_{M, N}}"]
			\arrow[dr, "{\gamma_{M, N}}"']
			& \Hom_{V \otimes V} \ab\big(M \otimes N, \;   \mathcal{R}_{\A, \A} (W^\prime))
			\arrow[d, "{\Hom_{V\otimes V}(M \otimes N,  \; \ev_{\xi})}"]
			\\
			& \Hom_{V \otimes V} \ab\big(M\otimes N, \;\mathcal{S}_{\A, \A} (\mathfrak{D}_{\xi}(W)))
		\end{tikzcd}
		\]
		Using the definition of $\psi^{W^\prime}_{M, N}$ in Theorem~\ref{thm:hom-realization}, the definition of $\gamma_{M, N}$ in Corollary~\ref{cor:li-intertwiner-hom},
		and the construction of a $V\otimes V$-module homomorphism from a natural transformation given in
		Lemma~\ref{lem:external-tensor-density}, we get
		\[
		\ev_{\xi} \ab\big(\Y(a, y) b) = \Y_{\xi}(a \otimes b),
		\]
		where $\Y$ is an intertwining operator of type $\binom{W^\prime}{M \; N}$
		for some generalized $V$-modules $M$ and $N$, with $a \in M$ and $b \in N$.
		This completes the proof.
		
	\end{proof}
	
	\begin{corollary}
		Let $V$ be a $C_2$-cofinite vertex operator algebra with nonnegative integer weights, and $W$ be a restricted $V$-module.
		Then there is a canonical isomorphism of restricted $V \otimes V$-modules
		\[
		\mathcal{R}(W^\prime) \cong \mathfrak{D}_{\xi}(W).
		\]
		
	\end{corollary}
	
	\begin{proof}
		By Proposition~\ref{prop:c2-regular-finiteness}, $\mathcal{R}(W^\prime)$ is a  restricted $V \otimes V$-module.
		Let $\mathcal{A}$ denote the category of generalized $V$-modules.
		Since $V$ is $C_2$-cofinite,  every weak $V \otimes V$-module is a bi-generalized $V \otimes V$-module.
		It follows from Proposition~\ref{prop:generated-category-characterizations}~(2) that
		\[
		\mathcal{R}_{\mathcal{A}, \mathcal{A}}(W^\prime) = \mathcal{R}(W^\prime), \quad \mathcal{S}_{\A, \A} \ab\big(\mathfrak{D}_{\xi}(W)) \cong \mathfrak{D}_{\xi}(W).
		\]
		The result now follows from Theorem~\ref{thm:li-regular-comparison}.
	\end{proof}

	\section*{Acknowledgement}
	
	Chao Yang is supported by the National Natural Science Foundation of China (Grant No. 12301039) and the Sichuan Science and Technology Program (Grant No. 2025ZNSFSC0797).


\end{document}